\documentclass[12pt]{article}
\usepackage{amssymb,amsmath}
\usepackage[all]{xy}

\usepackage[in]{fullpage}
\usepackage{sectsty}
\sectionfont{\large}
\subsectionfont{\normalsize}

\usepackage{amsthm}

\theoremstyle{plain}
\newtheorem{theorem}{Theorem}[section]
\newtheorem{lemma}[theorem]{Lemma}
\newtheorem{corollary}[theorem]{Corollary}
\newtheorem{proposition}[theorem]{Proposition}

\theoremstyle{definition}
\newtheorem{definition}[theorem]{Definition}
\newtheorem{example}[theorem]{Example}

\theoremstyle{remark}
\newtheorem{remark}[theorem]{Remark}

\newcommand{\lp}{\left(} 
\newcommand{\rp}{\right)} 
\newcommand{\lb}{\left\{} 
\newcommand{\rb}{\right\}} 
\newcommand{\la}{\left\langle } 
\newcommand{\ra}{\right\rangle} 
\newcommand{\ls}{\left[} 
\newcommand{\rs}{\right]} 

\newcommand{\italic}[1]{{\it #1\/}}
\newcommand{\dash}{\mathchar `\-}

\newcommand{\holders}{H\"{o}lder's }

\newcommand{\figa}{Fig\`a}

\newcommand{\tnorm}[1]{
  \left\vert\kern-0.9pt\left\vert\kern-0.9pt\left\vert #1
    \right\vert\kern-0.9pt\right\vert\kern-0.9pt\right\vert}

\newcommand{\norm}[1]{\left\| #1 \right\|}
\renewcommand{\mod}[1]{\left| #1 \right|}

\newcommand{\tensor}{\,\widehat{\otimes}\,}

\newcommand{\im}{\mathop{\rm Im}}

\newcommand{\dd}{\,{\rm d}}

\newcommand{\inj}{{\rm inj}}

\def\lin{\mathop{\rm lin\,}}

\def\im{\mathop{\rm im\,}}
\def\ker{\mathop{\rm ker\,}}

\def\sup{\mathop{\rm sup\,}}
\def\supp{\mathop{\rm supp\,}}

\def\#{\flat}

\newcommand {\C}{\mathbb{C}}

\newcommand {\R}{\mathbb{R}}

\newcommand {\N}{{\mathbb N}}
\newcommand {\F}{\mathbb{F}}

\newcommand {\B}{{\mathcal B}}

\newcommand {\AB}{{}_A\mathcal{B}}

\newcommand {\Amod}{A\dash{\bf mod}}

\newcommand {\Lmod}{L^{1}(G)\dash{\bf mod}}
\newcommand {\lmod}{\ell^{\, 1}(G)\dash{\bf mod}}

\newif\ifmargin 
\marginfalse 

\begin{document}
\begin{center}
{\Large Multi-norms and modules over group algebras} \\ \vspace{0.5cm} Paul Ramsden
\end{center}

\begin{abstract}
\noindent Let $G$ be a locally compact group, and let $1< p < \infty$. In this paper we investigate the injectivity of the left $L^1(G)$-module $L^p(G)$. We define a family of amenability type conditions called {\it $(p,q)$-amenability}, for any $1\leq p\leq q$. For a general locally compact group $G$ we show if $L^p(G)$ is injective, then $G$ must be $(p,p)$-amenable. For a discrete group $G$ we prove that $\ell^{\, p}(G)$ is injective if and only if $G$ is $(p,p)$-amenable.
\end{abstract}

\section{Introduction}
Let $G$ be a locally compact group, and let $L^{1}(G)$ be the group algebra of $G$. In \cite{DP}, H. G. Dales and M. E. Polyakov investigated when various canonical modules over $L^{1}(G)$ have certain well-known homological properties. One of the most difficult questions they considered seems to be to characterize the locally compact groups $G$ such that the left $L^{1}(G)$-module $L^{p}(G)$ is injective (for $1<p< \infty$). Since $L^p(G)$ is a dual Banach $L^1(G)$-module, it follows from Johnson's theorem that $L^{p}(G)$ is an injective left $L^{1}(G)$-module whenever $G$ is amenable. In \cite{DP} the authors obtained a partial converse in the case where $G$ is discrete. They showed in \cite[Theorem 5.12]{DP} that, if $G$ is a discrete group and $\ell^{\, p}(G)$ is injective for some $p \in (1,\infty)$, then $G$ must be `pseudo-amenable', a property very close to amenability.

In subsequent work Dales and Polyakov introduced the concept of a multi-normed space, and used this to define another generalized notion of amenability, now called {\it $(1,p)$-amenability}. They showed that, if $\ell^{\, p}(G)$ is injective, then $G$ must be $(1,p)$-amenable.

In this paper we define another generalized notion of amenability, called {\it $(p,q)$-amenability}, for any $1\leq p\leq q$. We prove for a discrete group $G$, that $\ell^{\, p}(G)$ is injective if and only if $G$ is {\it $(p,p)$-amenable}. For a general locally compact group $G$ we show the following implications
$$
L^{p}(G) \text{~injective} \Rightarrow G ~(p,p)\dash\text{amenable} \Rightarrow G ~(1,p)\dash\text{amenable~} \Rightarrow G \text{~pseudo}\dash\text{amenable}\,.
$$

\section{Preliminaries}
\subsubsection*{Banach spaces}
For $n \in \N$, we set $\N_n=\{1,\ldots,n\}$. The \italic{indicator function} of a subset $T$ of a set $S$ is denoted by $\chi_{T}$. We set $\delta_{s}=\chi_{\{s\}}\,\; (s \in S)$. Let $1\leq p\leq \infty$. The conjugate to $p$ is denoted by $p'$, so that $1\leq p'\leq \infty$ and satisfies $1/p+1/p'=1$.

Let $E$ be a Banach space; the convex hull of a subset $X\subset E$ is denoted by $\la X\ra$; the unit ball of $E$ is denoted by $E_{[1]}$; the identity operator on $E$ is denoted by $I_E$. We denote the dual space by $E'$; the action of $\lambda \in E'$ on an element $x \in E$ is written as $\la x, \lambda\ra$. Let $E$ and $F$ be Banach spaces. Then $\B(E, F)$ is the Banach space of all bounded linear operators from $E$ to $F$; the adjoint of an operator $T\in \B(E, F)$ is $T'\in \B(F', E')$. For $\lambda \in E'$ and $x \in F$ we define the rank $1$ operator $x \otimes \lambda\in \B(E, F)$ by
$$
(x \otimes \lambda)(y)=\la y , \lambda\ra x\quad (y \in E)\,.
$$

For a vector space $E$ and $n \in \N$, we denote by $E^n$ the vector space direct sum of $n$ copies of $E$. Following the notation of \cite{JAM} we define the \italic{weak $p$-summing norm} (for $1\leq p<\infty$) on $E^n$ by
$$
\mu_{p,n}(x)=\sup \left\{\left(\sum_{i=1}^{n}\mod{\left\langle x_{i}, \lambda\right\rangle}^{p}\right)^{1/p}: \lambda \in E'_{[1]}\right\}\index{$\mu_{p,n}(x)$}
$$
where $x=(x_1,\ldots,x_n) \in E^n$. We set $\ell^{\,p}_n(E)^{w}=(E^n, \mu_{p,n})$.

Let $K$ be a compact Hausdorff space. Then $C(K)$ is the Banach space of complex-valued continuous functions on $K$ equipped with the uniform norm $\mod{\,\cdot\,}_{K}$, given by
$$
\mod{f}_{K}=\sup\lb \mod{f(x)}: x \in \Omega\rb\quad (f \in C(K))\,.
$$

Let $\Omega$ be a locally compact space. Then $C_{00}(\Omega)$ is the normed space of continuous functions on $\Omega$ with compact support, again equipped with the uniform norm.

Again, let $\Omega$ be a locally compact space. We consider $\Omega$ as a measurable space equipped with the Borel $\sigma$-algebra on $\Omega$, denoted by $\mathcal{B}_{\Omega}$. Then $M(\Omega)$ is the Banach space of complex-valued regular Borel measures on $\Omega$, equipped with the total variation norm $\norm{\,\cdot\,}$, given by
$$
\norm{\mu}=\mod{\mu}(\Omega)\quad (\mu \in M(\Omega))\,.
$$

Let $\mu$ be a positive measure on $\Omega$, and let $1\leq p<\infty$. Then $L^p(\Omega)=L^p(\Omega, \mu)$ is the Banach space of complex-valued $p$-integrable functions on $\Omega$, equipped with the norm $\norm{\,\cdot\,}_p$, given by
$$
\norm{f}_p=\lp \int \mod{f}^p \dd \mu\rp^{1/p}\quad (f \in L^p(\Omega))\,.
$$

\subsubsection*{Banach homology}
For the homological background we refer the reader to the standard reference \cite{HE86}. We briefly sketch what we need. Let $A$ be a Banach algebra, and denote by $\Amod$ the category of Banach left $A$-modules. A module $E \in \Amod$ is \italic{faithful} if $x=0$ whenever $a\cdot x=0$ for all $a \in A$. For $E, F \in \Amod$ the space of $A$-module morphisms from $E$ to $F$ is $\AB(E, F)$. A monomorphism $T \in \AB(E, F)$ is \italic{admissible} if there exists $S \in {\cal B}(F,E)$ with \mbox{$S\,\circ\,T =I_E$}, and $T$ is a {\it coretraction} if there exists $S \in \AB(F, E)$ such that $S\circ T=I_{E}$.

\begin{definition} \label{1.1b}
Let $A$ be a Banach algebra, and let $J\in \Amod$.  Then $J$ is {\it injective} if, for each $E,F\in \Amod$, for each admissible monomorphism $T\in \AB(E,F)$, and for each $S\in \AB(E,J)$, there exists $R\in \AB(F,J)$ with $R\,\circ\,T=S$.
\end{definition}\smallskip

Let $A$ be a Banach algebra, and let $E$ be a Banach space. Then ${\cal B}(A,E) \in \Amod$ with the module operation
$$
(a\cdot T)(b)=T(ba)\quad (a, b\in A,\,T\in \B(A, E))\,.
$$
Now let $E \in \Amod$ be faithful. We define the \italic{canonical embedding} $\Pi:E\rightarrow {\cal B}(A,E)$ by the formula
$$
\Pi(x)(a)=a\cdot x\quad (a\in A,\,x\in E)\,.
$$

We shall use the following characterization of injective modules.

\begin{proposition}[{\cite[Proposition 1.7]{DP}}]\label{inj} \ifmargin \marginpar{inj} \fi
Let $A$ be a Banach algebra, and let $E\in \Amod$ be faithful. Then $E$ is injective if and only if the morphism $\Pi\in{\AB}(E, \B(A, E))$ is a coretraction.\qed
\end{proposition}

\subsection*{Multi-normed spaces}
The following definition is due to Dales and Polyakov. For a full account of the theory of multi-normed spaces see \cite{DP08}.

\begin{definition}[{\cite[Definition 2.1]{DP08}}]
Let $(E, \norm{\,\cdot\,})$ be a Banach space, and let $(\norm{\,\cdot\,}_n: n \in \N)$ be a sequence such that $\norm{\,\cdot\,}_n$ is a norm on $E^n$ for each $n \in \N$, with $\norm{\,\cdot\,}_1=\norm{\,\cdot\,}$ on $E=E^1$. Then the sequence $(\norm{\,\cdot\,}_n: n \in \N)$ is a \italic{multi-norm} if the following axioms hold (where in each case the axiom is required to hold for all $n\geq 2$ and all $x_1,\ldots, x_n \in E$):\smallskip

\noindent(A1) $\norm{(x_{\sigma(1)},\ldots,x_{\sigma(n)})}_{n}=\norm{(x_{1},\ldots,x_{n})}_{n}$ for each permutation $\sigma$;\smallskip

\noindent(A2) $\norm{(\alpha_1x_{1},\ldots,\alpha_nx_{n})}_{n}\leq \max_{i \in \N_n}\mod{\alpha_i}\norm{(x_{1},\ldots,x_{n})}_{n} \quad (\alpha_1,\ldots,\alpha_n \in \C^n)$;\smallskip

\noindent(A3) $\norm{(x_1,\ldots,x_{n-1},0)}_{n}=\norm{(x_1,\ldots,x_{n-1})}_{n-1}$;

\noindent(A4) $\norm{(x_1,\ldots,x_{n-2},x_{n-1},x_{n-1})}_{n}=\norm{(x_1,\ldots,x_{n-2},x_{n-1})}_{n-1}$.

The Banach space $E$ equipped with a multi-norm is a \italic{multi-normed space}.
\end{definition}

Suppose that in the above definition we replace axiom (A4) by the following axiom:\smallskip

\noindent(B4) $\norm{(x_1,\ldots,x_{n-2},x_{n-1},x_{n-1})}_{n}=\norm{(x_1,\ldots,x_{n-2},2x_{n-1})}_{n-1}$.\smallskip

\noindent Then we obtain the definition of a {\it dual multi-norm} and of a {\it dual multi-normed space}.

The following Lemma, whose proof we omit, is an elementary consequence of the axioms.

\begin{lemma}[{\cite[Lemma 4.2.8]{PRthesis}}]\label{4.2.8}
Let $E$ be a Banach space equipped with a sequence of norms $(\norm{\,\cdot\,}_n: n \in \N)$ on the spaces $E^n\,\; (n \in \N)$ which satisfy axioms {\rm(A1)-(A3)}. Let $n \geq 2$, and let $x_1,\ldots, x_{n-1} \in E$. Then:\smallskip

{\rm (i)} $\norm{(x_1,\ldots,x_{n-2}, x_{n-1}, x_{n-1})}_{n}\geq\norm{(x_1,\ldots,x_{n-2}, x_{n-1})}_{n-1}$;\smallskip

{\rm (ii)} $\norm{(x_1,\ldots,x_{n-2}, x_{n-1}, x_{n-1})}_{n} \leq\norm{(x_1,\ldots,x_{n-2}, 2x_{n-1})}_{n-1}$.\qed
\end{lemma}

On each Banach space $E$ there exists a unique {\it maximum multi-norm} $(\norm{\,\cdot\,}_n^{\max}: n \in \N)$ with the property that $\norm{x}_n\leq \norm{x}_n^{\max}\,\; (n \in \N,\, x \in E^n)$ for every multi-norm $(\norm{\,\cdot\,}_n: n \in \N)$ over $E$. By \cite[Proposition 3.23]{DP08}, for each $n \in \N$ and $x=(x_1,\ldots,x_n) \in E^n$, we have
\begin{equation}\label{max-mn}
\norm{x}^{\max}_n=\sup \lb \mod{\sum_{i=1}^n\la x_i, \lambda_i\ra}: \lambda_1,\ldots, \lambda_n \in E',\, \mu_{1,n}(\lambda_1,\ldots, \lambda_n)\leq 1\rb\,.
\end{equation}

Let $E$ and $F$ be Banach spaces, let $T:E\rightarrow F$ be a linear mapping, and let $k \in \N$. Then we define the $k^{\rm th}$-\italic{amplification} of $T$, $T^{(k)}:E^k\rightarrow F^k$ by\index{$T^{(k)}$}
$$
T^{(k)}(x_{1},\ldots,x_k)=(T(x_1),\ldots,T(x_k))\quad ((x_1,\ldots, x_k) \in E^k)\,.
$$
\begin{definition}
Let $E$ and $F$ be multi-normed spaces, and let $T \in \mathcal{B}(E, F)$. Then $T$ is \italic{multi-bounded} if 
$$
\norm{T}_{mb}=\sup_{k \in \N} \norm{T^{(k)}}_{\B(E^k, F^k)}< \infty\,.\index{$T_{mb}$}
$$
We set $\mathcal{M}(E, F)=\{T \in \B(E, F): \norm{T}_{mb}<\infty\}$. Then $\norm{\,\cdot\,}_{mb}$ is a norm on $\mathcal{M}(E, F)$ called the \italic{multi-bounded norm} and $\mathcal{M}(E, F)$ is the space of {\it multi-bounded operators}.\index{$M(E, F)$}
\end{definition}

It is easy to check that $\mathcal{M}(E, F)$ is a Banach space. This definition of the multi-bounded norm agrees with that given in \cite[Definition 5.11]{DP08} in terms of multi-bounded sets. In the case where $E$ and $F$ are operator sequence spaces, the multi-bounded norm is the same as the \italic{sequentially-bounded norm}. In this case, the multi-bounded operators are the same as the \italic{sequentially bounded operators}. See \cite{LNR04} for more about operator sequence spaces.

\begin{definition}[{\cite[Definition 6.3]{DP08}}]
Let $E$ be a multi-normed space. A subset $B\subset E$ is {\it multi-bounded} if
$$
c_B:=\sup \lb \norm{(x_1,\ldots, x_n)}_n: x_1,\ldots, x_n \in B,\, n \in \N\rb<\infty\,.
$$
\end{definition}

Let $E, F$ be multi-normed spaces, and let $T \in \mathcal{M}(E, F)$. It follows immediately from the definitions that $T(B)$ is a multi-bounded set in $F$ whenever $B$ is a multi-bounded set in $E$. Conversely, it is proved in \cite[Proposition 6.10]{DP08} that any $T \in \B(E, F)$ which takes multi-bounded sets to multi-bounded sets is multi-bounded, and further
$$
\norm{T}_{mb}=\sup\lb c_{T(B)}: c_B\leq 1\rb\,.
$$

\section{The weak $(p,q)$-multi-norm}
We now introduce the main class of examples of multi-norms. Let $E$ be a Banach space, and let $1\leq p,q< \infty$. For each $n \in \N$ and $x=(x_1,\ldots, x_n) \in E^n$ we set
$$
\norm{x}^{(p,q)}_n=\sup \lb \lp\sum_{i=1}^n\mod{\la x_i, \lambda_i\ra}^{q}\rp^{1/q}: \lambda \in (E')^n,\, \mu_{p,n}(\lambda)\leq 1\rb\,.
$$
It is clear that $\norm{\,\cdot\,}_n$ is a norm on $E^n$. 

\begin{proposition}\label{4.6.4}
Let $E$ be a Banach space, and let $1\leq p\leq q< \infty$. Then the family $(\norm{\,\cdot\,}_n^{(p,q)}: n\in \N)$ is a multi-norm on $E$.
\end{proposition}
\begin{proof}
It is clear that the family $(\norm{\,\cdot\,}_n^{(p,q)}: n\in \N)$ satisfies axioms (A1)-(A3). We shall verify axiom (A4).

Let $n\geq 2$, let $x_1,\ldots, x_{n-1} \in E$, and let $\varepsilon >0$. Set $x=(x_1,\ldots, x_{n-2},x_{n-1},x_{n-1}) \in E^n$. There exists $\lambda_1,\ldots, \lambda_n \in E'$ with $\mu_{n,p}(\lambda_1,\ldots,\lambda_n)\leq 1$, and such that
$$
\norm{x}^{(p,q)}_n-\varepsilon<\lp\sum_{i=1}^{n-2}\mod{\la x_i, \lambda_i\ra}^{q}+\mod{\la x_{n-1}, \lambda_{n-1}\ra}^{q}+\mod{\la x_{n-1}, \lambda_{n}\ra}^{q}\rp^{1/q}\,.
$$

By the Hahn--Banach theorem there exists $\alpha, \beta \in \C$ with $\norm{(\alpha, \beta)}_{q'}\leq 1$ and
$$
\mod{\la x_{n-1}, \lambda_{n-1}\ra}^{q}+\mod{\la x_{n-1}, \lambda_{n-1}\ra}^{q}=\la x_{n-1}, \alpha\lambda_{n-1}+\beta\lambda_n\ra^{q}\,.$$

Set $\gamma=\norm{(\alpha, \beta)}_{p'}$. Since $q'\leq p'$ we have $\gamma\leq \norm{(\alpha, \beta)}_{q'}\leq 1$. We have
\begin{equation*}
\begin{split}
\mu_{p,n-1}(\lambda_1,\ldots, \lambda_{n-2},\alpha\lambda_{n-1}+\beta\lambda_n)&\leq \mu_{p,n}(\lambda_1,\ldots, \lambda_{n-2},\gamma\lambda_{n-1},\gamma\lambda_n) \\
&\leq \max\{1, \gamma\}\mu_{p,n}(\lambda_1,\ldots,\lambda_n)\leq 1\,.
\end{split}
\end{equation*}
Hence
$$
\norm{x}^{(p,q)}_n-\varepsilon<\lp\sum_{i=1}^{n-2}\mod{\la x_i, \lambda_i\ra}^{q}+\la x_{n-1}, \alpha\lambda_{n-1}+\beta\lambda_n\ra^{q}\rp^{1/q}\leq \norm{(x_1,\ldots,x_{n-1})}^{(p,q)}_{n-1}\,.
$$
The reverse inequality is given by Lemma \ref{4.2.8}(i).
\end{proof}

\begin{definition}
Let $E$ be a Banach space, and let $1\leq p\leq q< \infty$. Then the \italic{weak $(p,q)$-multi-norm} over $E$ is the multi-norm $(\norm{\,\cdot\,}^{(p,q)}_n: n \in \N)$ described above.
\end{definition}

\begin{remark}
The weak-$(1,1)$ multi-norm is just the maximum multi-norm.
\end{remark}

The next result is a straight forward verification from the definitions.

\begin{proposition}\label{4.6.7}
Let $1\leq p\leq q<\infty$, and let $E, F$ be multi-normed spaces both equipped with the weak $(p,q)$-multi-norm. Then $\mathcal{M}(E, F)=\B(E, F)$ with $\norm{T}_{mb}=\norm{T}\,\; (T \in \mathcal{M}(E, F))$.\qed
\end{proposition}

\begin{lemma}\label{4.6.6}\ifmargin \marginpar[4.6.6]{4.6.6} \fi
Let $E$ be a Banach space, and let $1\leq p,q< \infty$. Then for each $n \in \N$ and $\lambda\in (E')^n$ we have
\begin{multline*}
\sup \lb \lp\sum_{i=1}^n\mod{\la \lambda_i, \Phi_i\ra}^{q} \rp^{1/q}: \Phi \in (E'')^n,\, \mu_{p,n}(\Phi)\leq 1\rb \\
=\sup \lb \lp\sum_{i=1}^n\mod{\la x_i, \lambda_i\ra}^{q} \rp^{1/q}: x \in E^n,\, \mu_{p,n}(x)\leq 1\rb\,.
\end{multline*}
\end{lemma}
\begin{proof}
We denote the left and right hand sides of this equation by $\norm{\lambda}$ and $\tnorm{\lambda}$ respectively. Clearly $\norm{\lambda}\geq\tnorm{\lambda}$. Take $\varepsilon>0$ and $\Phi=(\Phi_1,\ldots, \Phi_n) \in (E'')^n$ with $\mu_{p,n}(\Phi)\leq 1$ and
$$
\norm{\lambda}-\varepsilon\leq \lp\sum_{i=1}^n\mod{\la \lambda_i, \Phi_i\ra}^{q}\rp^{1/q}\,.
$$
By the Principle of Local Reflexivity there exists $T \in \B(E'', E)$ with $\norm{T}\leq 1+\varepsilon$ and
$$
\la T(\Phi_i), \lambda_i\ra=\la \lambda_i, \Phi_i\ra\quad (i \in \N_n)\,.
$$
Set $x=(T(\Phi_1),\ldots, T(\Phi_n)) \in E^n$. We have
$$
\mu_{p,n}(x)\leq(1+\varepsilon)\mu_{p,n}(\Phi)\leq 1+\varepsilon\,.
$$
Now we have
$$
\tnorm{\lambda}\geq \frac{1}{1+\varepsilon}\lp\sum_{i=1}^n\mod{\la T(\Phi_i), \lambda_i\ra}^{q}\rp^{1/q}=\frac{1}{1+\varepsilon}\lp\sum_{i=1}^n\mod{\la \lambda_i, \Phi_i\ra}^{q}\rp^{1/q}\geq \frac{1}{1+\varepsilon}\lp\norm{\lambda}-\varepsilon\rp\,.
$$
This is true for each $\varepsilon>0$. Therefore $\tnorm{\lambda}\geq \norm{\lambda}$ as required.
\end{proof}

\subsection{The dual of the weak $(p,q)$-multi-norm}
In this section we study the dual of the the norm $\norm{\,\cdot\,}_n^{(p,q)}$ on $(E')^n$.

Let $E$ be a Banach space, and let $n \in \N$. For $\alpha=(\alpha_i) \in \C^n$ we define the operator $M_{\alpha}:E^n\rightarrow E^n$ by
$$
M_{\alpha}(x_1,\ldots, x_n)=(\alpha_1x_{1},\ldots,\alpha_nx_{n})\,.
$$
Let $1\leq r< \infty$ and $1< s\leq \infty$. For each $n \in \N$ and $x \in E^n$ we set
$$
\tnorm{x}^{(r,s)}_n=\inf\lb \sum_{k=1}^N\norm{\alpha_{k}}_{s}\mu_{r,n}(y_{k})\rb\,,
$$
where the infimum is taken over all representations
$$
x=\sum_{k=1}^N M_{\alpha_k}(y_k)
$$
with $\alpha_k \in \C^n,\, y_k \in E^n,\, k \in \N_N,\, N \in \N$. It is clear that $\tnorm{\,\cdot\,}_n^{(r,s)}$ is a norm on $E^n$.

\begin{proposition}
Let $E$ be a Banach space, and let $1\leq r, s\leq \infty$ with $1<s\leq r'\leq \infty$. Then the family $(\tnorm{\,\cdot\,}_n^{(r,s)}: n\in \N)$ is a dual multi-norm over $E$.
\end{proposition}
\begin{proof}
This is `dual' to the proof of Proposition \ref{4.6.4}. We shall verify axiom (B4). Take $n\geq 2$, $x_1,\ldots,x_{n-1} \in E$ and a representation
$$
(x_1,\ldots,x_{n-2},x_{n-1},x_{n-1})=\sum_{k=1}^N M_{\alpha_k}(y_k)
$$
with $\alpha_k \in \C^n,\, y_k \in E^n,\, k \in \N_N,\, N \in \N$. For each $k \in \N_N$ we set $$\gamma_{k}=\norm{(\alpha_{k,n-1},\alpha_{k,n})}_s\,,$$
$$
z_k=\lp y_{k,1},\ldots,y_{k,n-2},\frac{\alpha_{k,n-1}y_{k,n-1}+\alpha_{k,n}y_{k,n}}{\gamma_k}\rp \in E^{n-1}\,,
$$
and
$$
\beta_k=(\alpha_{k,1},\ldots,\alpha_{k,n-2},\gamma_{k})\in \C^{n-1}\,.
$$

Then we have
$$
(x_1,\ldots,x_{n-2},2x_{n-1})=\sum_{k=1}^N M_{\beta_k}(z_k)\,.
$$
Further, we have $\norm{\beta_k}_{s}=\norm{\alpha_k}_{s}\,\; (k \in \N)$, and
\begin{equation*}
\begin{split}
\mu_{r,n-1}(z_{k})&\leq \max\lb 1, \frac{\norm{(\alpha_{k,n-1},\alpha_{k,n})}_{r'}}{\gamma_k}\rb\mu_{r,n}(y_k) \\
&=\mu_{r,n}(y_k)\quad \text{(since } s\leq r'\text{)}\,.
\end{split}
\end{equation*}

Hence $$\tnorm{(x_1,\ldots,x_{n-2},2x_{n-1})}_{n-1}^{(r,s)}\leq \sum_{k=1}^N\norm{\alpha_{k}}_{r}\mu_{r,n}(y_{k})\,.$$
Since this holds for all representations of $(x_1,\ldots,x_{n-2},x_{n-1},x_{n-1})$ we have
$$
\tnorm{(x_1,\ldots,x_{n-2},2x_{n-1})}_{n-1}^{(r,s)}\leq\tnorm{(x_1,\ldots,x_{n-2},x_{n-1},x_{n-1})}_{n}^{(r,s)}\,.
$$
The reverse inequality is given by Lemma \ref{4.2.8}(ii).
\end{proof}

\begin{definition}
Let $E$ be a Banach space, and let $1\leq r, s\leq \infty$ with $1<s\leq r'\leq \infty$. Then the \italic{weak $(r,s)$-dual multi-norm} over $E$ is the dual multi-norm $(\tnorm{\,\cdot\,}^{(r,s)}_n: n \in \N)$ described above.
\end{definition}

Let $E$ be a Banach space, and take $1\leq p,q < \infty$. Define an embedding $\nu_{E}:E^n\rightarrow \B(\ell^{\,p}_n(E')^{w}, \ell^{\, q}_n)$ by
$$
\nu_{E}(x)(\lambda)=\lp\la x_1, \lambda_1\ra,\ldots,\la x_n, \lambda_n\ra\rp\,,
$$
where $x=(x_1,\ldots,x_n) \in E^n$ and $\lambda=(\lambda_1,\ldots,\lambda_n) \in \ell^{\,p}_n(E')^{w}$ . Then there is an isometric embedding
$$
\nu_{E}:\lp E^n, \norm{\,\cdot\,}_n^{(p,q)}\rp\rightarrow \B(\ell^{\,p}_n(E')^{w}, \ell^{\, q}_n)\,.
$$

Now take $1\leq r< \infty$ and $1< s\leq \infty$. Define a map $\theta_{E}:\ell^{\,r}_n(E)^{w}\tensor \ell^{\, s}_n\rightarrow E^n$ by
$$
\theta_E(x\otimes \alpha)=M_{\alpha}(x)=(\alpha_1x_1,\ldots,\alpha_nx_n)\,,
$$
where $x=(x_1,\ldots,x_n) \in E^n$ and $\alpha=(\alpha_1,\ldots, \alpha_n)\in \ell^{\, s}_n$. The map $\theta_E$ is a linear surjection with closed kernel. There is an isometric isomorphism of Banach spaces
$$
\ell^{\,r}_n(E)^{w}\tensor \ell^{\, s}_n/\ker \theta_E\cong \lp E^n, \tnorm{\,\cdot\,}_n^{(r,s)}\rp\,.
$$

\begin{theorem}\label{duality}
Let $E$ be a Banach space, let $1\leq p, q, r< \infty$, and let $1< s\leq \infty$. Then there are isometric isomorphisms of Banach spaces:\smallskip
 
{\rm (i)} $\lp E^n, \norm{\,\cdot\,}_n^{(p,q)}\rp'\cong \lp(E')^n, \tnorm{\,\cdot\,}_n^{(p,q')}\rp$; and\smallskip

{\rm (ii)} $\lp E^n, \tnorm{\,\cdot\,}_n^{(r,s)}\rp'\cong \lp(E')^n, \norm{\,\cdot\,}_n^{(r,s')}\rp$.
\end{theorem}
\begin{proof}
(i) It is easily checked that the following diagram commutes
\begin{equation*}
\SelectTips{eu}{12}\xymatrix{(\ell^{\,p}_n(E')^{w}\tensor\ell^{\, q'}_n)'' \ar[r]^-{(\nu_{E})'} & (E')^n  \\ 
                             \ell^{\,p}_n(E')^{w}\tensor\ell^{\, q'}_n\,.\ar[u] \ar[ur]_{\theta_{E'}}& ~}
\end{equation*}
Hence we have isometric isomorphisms of Banach spaces
\begin{equation*}
\begin{split}
\lp E^n, \norm{\,\cdot\,}_n^{(p,q)}\rp'&\cong (\ell^{\,p}_n(E')^{w}\tensor\ell^{\, q'}_n)''/\ker (\nu_{E})' \\
&\cong \ell^{\,p}_n(E')^{w}\tensor \ell^{\, q'}_n/\ker \theta_{E'}\cong \lp(E')^n, \tnorm{\,\cdot\,}_n^{(p,q')}\rp\,.
\end{split}
\end{equation*}

(ii) Similarly, the following diagram commutes
\begin{equation*}
\SelectTips{eu}{12}\xymatrix{(E')^n \ar[r]^-{(\theta_{E})'} \ar[dr]_{\nu_{E'}} & \B(\ell^{\,r}_n(E)^{w}, \ell^{\, s'}_n)  \\ 
            ~ & \B(\ell^{\,r}_n(E'')^{w}, \ell^{\, s'}_n)\,. \ar[u]_{j: T\mapsto T|E^n}}
\end{equation*}
Hence there is an isometric isomorphism
$$
\lp E^n, \tnorm{\,\cdot\,}_n^{(r,s)}\rp'\cong \im (\theta_{E})'=\im j\circ \nu_{E'}\,.
$$
By Lemma \ref{4.6.6}, there is an isometric isomorphism 
\begin{equation*}
\im j\circ \nu_{E'}\cong \im \nu_{E'}\cong \lp(E')^n, \norm{\,\cdot\,}_n^{(r,s')}\rp\,.\qedhere
\end{equation*}
\end{proof}

The next lemma follows easily from \holders inequality. 

\begin{lemma}\label{4.6.13a}
Let $1\leq p<\infty$, and let $E$ be a Banach space. For each $1< u, v<\infty$ with $1/u+1/v=1$, we have
\begin{equation*}
\mu_{p}(M_{\alpha}(x))\leq \norm{\alpha}_{pu}\mu_{pv}(x)\quad (x \in E^n,\, \alpha \in \C^n)\,.\vspace{-\baselineskip}
\end{equation*}\qed
\end{lemma}

Let $\alpha=(\alpha_1,\ldots,\alpha_n) \in \C^n$, and let $t \in \R\setminus \{0\}$. Then we set
$$
\alpha^t=(\alpha_1^t,\ldots,\alpha_n^t) \in \C^n\,.
$$

\begin{proposition}
Let $1\leq p_1\leq q_1<\infty$ and $1\leq p_2\leq q_2<\infty$, and let $E$ be a Banach space. Suppose that
$$
\text{{\rm (i)}}\,\; q_2\leq q_1\,,\quad\text{and}\quad {\rm (ii)}\,\; \frac{1}{p_2}-\frac{1}{q_2}\leq \frac{1}{p_1}-\frac{1}{q_1}\,.
$$
Then $\tnorm{\,\cdot\,}_n^{(p_2,q_2')}\leq \tnorm{\,\cdot\,}_n^{(p_1,q_1')}$ and $\norm{\,\cdot\,}_n^{(p_1,q_1)}\leq \norm{\,\cdot\,}_n^{(p_2,q_2)}$ on $E^n$.
\end{proposition}
\begin{proof}
We shall prove that $\tnorm{\,\cdot\,}_n^{(p_2,q_2')}\leq \tnorm{\,\cdot\,}_n^{(p_1,q_1')}$, the other inequality follows by duality.

{\bf Case 1} $p_1\leq p_2$. In this case the hypothesis (ii) is automatically true. The result follows from the inequalities
$$
\norm{\alpha}_{q_2'}\leq \norm{\alpha}_{q_1'}\,,\quad \mu_{p_2}(y)\leq \mu_{p_1}(y)\quad (\alpha \in \C^n,\, y \in E^n)\,.
$$

{\bf Case 2} $p_1>p_2$. Let $\alpha \in \C^n$ and $y \in E^n$. For each $t \in (0,1)$ and $u, v \in (1,\infty)$ with $1/u+1/v=1$, we have
\begin{equation*}
\begin{split}
\norm{\alpha^t}_{q_2'}\mu_{p_2}(M_{\alpha^{1-t}}(y))&=\norm{\alpha}^t_{tq_2'}\mu_{p_2}(M_{\alpha^{1-t}}(y)) \\
&\leq \norm{\alpha}^t_{tq_2'}\norm{\alpha^{1-t}}_{p_2u}\mu_{p_2v}(y)\quad\text{(by Lemma \ref{4.6.13a})} \\
&=\norm{\alpha}^t_{tq_2'}\norm{\alpha}_{(1-t)p_2u}^{1-t}\mu_{p_2v}(y)\,.
\end{split}
\end{equation*}

Suppose that $(t,u) \in (0,1)\times (1,\infty)$ can be chosen to satisfy the following inequalities
$$
tq_2'\geq q_1',\quad (1-t)p_2u\geq q_1',\quad p_2\lp\frac{u}{u-1}\rp\geq p_1\,.
$$
Then we have
$$\norm{\alpha^t}_{q_2'}\mu_{p_2}(M_{\alpha^{1-t}}(y))\leq \norm{\alpha}_{q_1'}\mu_{p_1}(y)\,.$$
Now take $x \in E^n$ and a representation $x=\sum_{k=1}^NM_{\alpha_k}(y_k)$ where $\alpha_k \in \C^n,\, y_k \in E^n\,\; (k \in \N_N)$. Then
$$
x=\sum_{k=1}^NM_{\alpha_k^{t}}\lp M_{\alpha_k^{1-t}}(y_k)\rp
$$ 
and
$$
\sum_{k=1}^N\norm{\alpha_k^{t}}_{q_2'}\mu_{p_2}\lp M_{\alpha_k^{1-t}}(y_k)\rp\leq \sum_{k=1}^N\norm{\alpha_k}_{q_1'}\mu_{p_1}(y_k)\,.
$$
It follows that $\tnorm{x}_n^{(p_2,q_2')}\leq \tnorm{x}_n^{(p_1,q_1')}$.

To complete the proof we need to show that such a choice of $(t, u)$ is possible. Indeed, set $t=\frac{q_1'}{q_2'}$ and $u=\frac{p_1}{p_1-p_2}$. Then $t \in (0,1)$ and $u>1$ by the hypothesis of case 2. For these choices we have
$$
tq_2'= q_1'\quad\text{and}\quad p_2\lp\frac{u}{u-1}\rp=p_1\,.
$$
The remaining inequality follows from a rearrangement of the inequality (ii). Indeed, (ii) is equivalent to
$$
\frac{1}{p_2}-\frac{1}{p_1}\leq \frac{1}{q_1'}-\frac{1}{q_2'}\,. 
$$
Multiplying by $p_2p_1q_1'$ gives
$$
(p_1-p_2)q_1'\leq \lp 1-\frac{q_1'}{q_2'}\rp p_2p_1=p_2p_1(1-t)\,.
$$
Since $p_1-p_2>0$ we can divide this out to get
$$
q_1'\leq\frac{p_2p_1(1-t)}{p_1-p_2}=p_2u(1-t)\,,
$$
as required.
\end{proof}

\begin{corollary}\label{4.6.15a}
Let $1<p< q<\infty$, and let $E$ be a Banach space. For each $n \in \N$ the following inequalities hold on $E^n$:\medskip

{\rm (i)} $\norm{\,\cdot\,}_n^{(1,q)}\leq \norm{\,\cdot\,}_n^{(p,q)}\leq \norm{\,\cdot\,}_n^{(q,q)}\,;$\medskip

{\rm (ii)} $\norm{\,\cdot\,}_n^{(q,q)}\leq \norm{\,\cdot\,}_n^{(p,p)}\leq \norm{\,\cdot\,}_n^{(1,1)}\,.$\qed
\end{corollary}

\subsection{The weak $(1,q)$-multi-norm over $L^{\, 1}(\Omega)$}
Let $\Omega$ be a locally compact measure space. In this section we give a concrete description of the weak $(1,q)$-multi-norm over the Banach spaces $L^1(\Omega)$ and $M(\Omega)$. We shall identify this norm with the \italic{standard $(1, q)$-multi-norm} defined in \cite[Definition 4.7]{DP08}. The result is based on the following identification of $\mu_{1,n}$.

\begin{proposition}[{\cite[2.6]{JAM}}]\label{4.6.13}
Let $K$ be a compact space, and let $\lambda_{1},\ldots,\lambda_{n} \in C(K)$. Then
\begin{equation*}
\mu_{1, n}(\lambda_1,\ldots, \lambda_n)=\mod{\sum_{i=1}^n\mod{\lambda_i}}_K\,.\vspace{-\baselineskip}
\end{equation*}\qed
\end{proposition}

Let $\Omega$ be a locally compact space. Then $M(\Omega)'$ is isometrically isomorphic to $C(\widetilde{\Omega})$ for some compact space $\widetilde{\Omega}$ called the \italic{hyper-stonean cover of $\Omega$}, (this follows from general C*-algebra theory as in \cite[III.2.3]{TA79}, a direct proof is given in \cite{DLS2}). There is an isometric embedding $\kappa:M(\Omega)\rightarrow M(\widetilde{\Omega})$ which identifies $M(\Omega)$ with the closed subspace of $M(\widetilde{\Omega})$ consisting of the \italic{normal} measures on $\widetilde{\Omega}$. Thus we shall interpret $\mu \in M(\Omega)$ as a measure on $\widetilde{\Omega}$. The duality between $M(\Omega)$ and $C(\widetilde{\Omega})$ is then given by
$$
\la \mu, \lambda\ra=\int_{\widetilde{\Omega}}\lambda \dd \mu\quad (\mu \in M(\Omega),\, \lambda \in C(\widetilde{\Omega}))\,.
$$
Let $\mu, \nu \in M(\Omega)$ be positive measures. Then we define $\mu\vee\nu \in M(\Omega)$ by
$$
(\mu\vee\nu)(E)=\sup\{\mu(E_1)+\nu(E_2): (E_1, E_2)\text{ a measurable partition of } E\}\,,
$$
where $E \in \B_{\Omega}$.

\begin{theorem}\label{4.6.14} \ifmargin \marginpar[4.6.14]{4.6.14} \fi
Let $\Omega$ be a locally compact space. Then the maximum multi-norm over $M(\Omega)$ is given by
$$
\norm{(\mu_1,\ldots, \mu_n)}_n^{\max}=\norm{\mod{\mu_1}\vee\cdots\vee\mod{\mu_n}}\quad (\mu_1,\ldots,\mu_n \in M(\Omega))\,.
$$
\end{theorem}
\begin{proof}
Take $n \in \N$ and $\mu_{1},\ldots,\mu_{n} \in M(\Omega)$. For $\lambda_{1},\ldots,\lambda_{n} \in C(\widetilde{\Omega})$ we have
\begin{equation*}
\begin{split}
\mod{\sum_{i=1}^{n}\la \mu_i, \lambda_i\ra}&\leq \sum_{i=1}^{n}\mod{\la \mu_i, \lambda_{i}\ra} \leq \sum_{i=1}^{n}\la \mod{\mu_i}, \mod{\lambda_{i}}\ra \leq \sum_{i=1}^{n}\la \mod{\mu_{1}}\vee\cdots\vee \mod{\mu_{n}}, \mod{\lambda_{i}}\ra \\
&=\la \mod{\mu_{1}}\vee\cdots\vee \mod{\mu_{n}}, \sum_{i=1}^n\mod{\lambda_i}\ra\leq\norm{\mod{\mu_{1}}\vee\cdots\vee \mod{\mu_{n}}}\mod{\sum_{i=1}^{n}\mod{\lambda_{i}}}_{\widetilde{\Omega}} \\
&=\norm{\mod{\mu_{1}}\vee\cdots\vee \mod{\mu_{n}}}\mu_{1, n}(\lambda_{1},\ldots,\lambda_{n})\,.
\end{split}
\end{equation*}
Hence
\begin{equation*}
\begin{split}
\norm{(\mu_1,\ldots, \mu_n)}_{n}^{\max}&=\sup\left\{\mod{\sum_{i=1}^{n}\left\langle \mu_i, \lambda_i\right\rangle}: \mu_{1, n}(\lambda_{1},\ldots,\lambda_{n})\leq 1\right\} \\
&\leq\norm{\mod{\mu_1}\vee\cdots\vee\mod{\mu_n}}\,.
\end{split}
\end{equation*}
Therefore $\norm{(\mu_1,\ldots, \mu_n)}_{n}^{\max}=\norm{\mod{\mu_1}\vee\cdots\vee\mod{\mu_n}}$, as required.
\end{proof}

Let $\Omega$ be a locally compact space, and let $X\subset \Omega$ be a measurable set. We define a projection $P_X:M(\Omega)\rightarrow M(\Omega)$ by
$$
P_X(\mu)(E)=\mu(X\cap E)\quad (E \in \B_{\Omega})\,.
$$
When restricted to $L^1(\Omega)$ this map has the form
$$
P_X(f)=\chi_{X}f\quad (f \in L^1(\Omega))\,.
$$

\begin{proposition}\label{4.6.15}
Let $\Omega$ be a locally compact space, and let $\mu_1,\ldots, \mu_n \in M(\Omega)$. Then:
\begin{equation*}
\norm{\mod{\mu_1}\vee\cdots\vee\mod{\mu_n}}=\sup_{\mathbf{X}}\sum_{i=1}^n\norm{P_{X_i}(\mu_i)}\,,
\end{equation*}
where the supremum is taken over all measurable partitions $\mathbf{X}=(X_1,\ldots,X_n)$ of $\Omega$.
\end{proposition}
\begin{proof}
Take $\mu_1,\ldots, \mu_n \in M(\Omega)$. For each measurable partition $\mathbf{X}=(X_1,\ldots,X_n)$ of $\Omega$, we have
$$\sum_{i=1}^n\norm{P_{X_i}(\mu_i)}=\sum_{i=1}^n\mod{\mu_i}(X_i)\leq \sum_{i=1}^n(\mod{\mu_1}\vee\cdots\vee \mod{\mu_{n}})(X_i)=\norm{\mod{\mu_1}\vee\cdots\vee \mod{\mu_{n}}}\,.$$
Hence $\sup_{\mathbf{X}}\sum_{i=1}^n\norm{P_{X_i}(\mu_i)}\leq \norm{\mod{\mu_1}\vee\cdots\vee \mod{\mu_{n}}}$.

We shall prove that for each $n \geq 2$ and $\mu_1,\ldots, \mu_n \in M(\Omega)$, there exists a measurable partition $(X_1,\ldots,X_n)$ of $\Omega$ with
$$
\norm{\mod{\mu_1}\vee\cdots\vee\mod{\mu_n}}=\sum_{i=1}^n\norm{P_{X_i}(\mu_i)}\,.
$$
This is most easily done by induction on $n$.

For positive measures $\mu, \nu \in M(\Omega)$, by the Hahn decomposition theorem (\cite[Theorem 4.1.4]{COHN}), there exists a set $(\mu \geq \nu) \in \B_{\Omega}$ with the property that $\mu(E)\geq\nu(E)$ for all measurable subsets $E\subset (\mu \geq \nu)$, and $\mu(E) \leq \nu(E)$ for all measurable subsets $E\subset \Omega\setminus (\mu \geq \nu)$.

Consider the case $n=2$, and let $\mu_1, \mu_2 \in M(\Omega)$. Set $X_1=(\mod{\mu_1}\geq \mod{\mu_2})$ and $X_2=\Omega\setminus X_1$. Then we have
\begin{equation*}
\begin{split}
\norm{\mod{\mu_1}\vee\mod{\mu_2}}&=\mod{\mu_1}\vee\mod{\mu_2}(X_1)+\mod{\mu_1}\vee\mod{\mu_2}(X_2) \\
&=\mod{\mu_1}(X_1)+\mod{\mu_2}(X_2) \\
&=\norm{P_{X_1}(\mu_1)}+\norm{P_{X_2}(\mu_2)}\,.
\end{split}
\end{equation*}

Now assume that the result holds for some $n \in \N$, and take $\mu_1,\ldots, \mu_{n+1} \in M(\Omega)$. Set $\mu=\mod{\mu_1}\vee\cdots\vee\mod{\mu_n}$ and $X=(\mu \geq \mod{\mu_{n+1}})$. By the inductive hypothesis there is a measurable partition $(Y_1,\ldots, Y_n)$ of $\Omega$ with
$$
\norm{\mod{P_{X}(\mu_1)}\vee\cdots\vee\mod{P_{X}(\mu_n)}}=\sum_{i=1}^n\norm{P_{Y_i}(P_{X}(\mu_i))}=\sum_{i=1}^n\norm{P_{Y_i\cap X}(\mu_i)}\,.
$$
Then we have
\begin{equation*}
\begin{split}
\norm{\mod{\mu_1}\vee\cdots\vee\mod{\mu_{n+1}}}&=\norm{\mu\vee\mod{\mu_{n+1}}}=\norm{P_{X}(\mu)}+\norm{P_{\Omega\setminus X}(\mu_{n+1})} \\
&=\norm{\mod{P_{X}(\mu_1)}\vee\cdots\vee\mod{P_{X}(\mu_n)}}+\norm{P_{\Omega\setminus X}(\mu_{n+1})} \\
&=\sum_{i=1}^n\norm{P_{Y_i\cap X}(\mu_i)}+\norm{P_{\Omega\setminus X}(\mu_{n+1})}\,,
\end{split}
\end{equation*}
where the sets $(Y_1\cap X,\ldots,Y_n\cap X,\Omega\setminus X)$ form a measurable partition of $\Omega$. By induction, the result follows.
\end{proof}

Let $K$ be a compact space, and take $n \in \N$. We define $D_n$ to be the set of $(\lambda_{1}, \ldots, \lambda_n) \in C(K)^n$ such that $\mod{\lambda_i}_K\leq 1\,\; (i \in \N_n)$, and the sets $\supp \lambda_1,\ldots, \supp \lambda_n$ are pairwise disjoint.

\begin{corollary}\label{4.6.16} \ifmargin \marginpar[4.6.15]{4.6.15} \fi
Let $K$ be a compact space. Then $(C(K)^n, \mu_{1,n})_{[1]}=\overline{\la D_n\ra}$, where the closure is in the weak-$*$ topology.
\end{corollary}
\begin{proof}
Set $B_n=(C(K)^n, \mu_{1,n})_{[1]}$. It is easily seen using Proposition \ref{4.6.13} that $\overline{\la D_n\ra}\subset B_n$. Assume towards a contradiction that there exists $\varphi \in B_n\setminus \overline{\la D_n\ra}$. Since $\overline{\la D_n\ra}$ is a balanced set, by a corollary to the Hahn--Banach separation theorem (\cite[Theorem 3.7]{RUDIN73}), there exists $\mu=(\mu_1,\ldots,\mu_n) \in M(K)^n$ with
$$
\la \lambda, \mu\ra \leq 1\quad (\lambda \in \overline{\la D_n\ra})\quad {\rm and}\quad \la \varphi, \mu\ra >1\,. 
$$
By Proposition \ref{4.6.15}, we have
$$
\norm{\mod{\mu_1}\vee\cdots\vee\mod{\mu_{n}}}=\sup_{\mathbf{X}}\sum_{i=1}^n\norm{P_{X_i}(\mu_i)}=\sup_{\lambda \in D_n}\mod{\sum_{i=1}^n\la \lambda_i,\mu_i\ra}\leq 1 < \norm{\mu}_{n}^{\max}\,,
$$
which is a contradiction of Theorem \ref{4.6.14}. Therefore $B_n=\overline{\la D_n\ra}$.
\end{proof}

\begin{theorem}\label{4.6.17} \ifmargin \marginpar[4.6.17]{4.6.17} \fi
Let $\Omega$ be a locally compact space, and let $1\leq q< \infty$. Then the weak $(1,q)$-multi-norm over $M(\Omega)$ is given by
\begin{equation*}
\norm{(\mu_1,\ldots,\mu_n)}_n^{(1,q)}=\sup_{\mathbf{X}}\lp \sum_{i=1}^n\norm{P_{X_i}(\mu_i)}^q\rp^{1/q} \quad (\mu_1,\ldots,\mu_n \in M(\Omega))\,.
\end{equation*}
where the supremum is taken over all measurable partitions $\mathbf{X}=(X_1,\ldots,X_n)$ of $\Omega$.
\end{theorem}
\begin{proof}
Take $\mu_1,\ldots ,\mu_n \in M(\Omega)$. By Corollary \ref{4.6.16} we have
\begin{multline*}
\sup \lb \lp \sum_{i=1}^n\mod{\la \mu_i, \lambda_i\ra}^q \rp^{1/q}: \lambda \in (C(\widetilde{\Omega})^n, \mu_{1,n})_{[1]} \rb= \\
\sup \lb \lp \sum_{i=1}^n\mod{\la \mu_i, \lambda_i\ra}^q \rp^{1/q}: \lambda \in  D_n\rb\,.
\end{multline*}
This gives the result.
\end{proof}

The following remark is contained in \cite[Example 4.9]{DP08}. Let $\Omega$ be a locally compact space. Then $L^1(\Omega)''$ is isometrically isomorphic to $M(K)$ for a certain compact space $K$. Let $X \in \B_{\Omega}$, and let $P_X \in \B(L^1(\Omega))$ be the projection onto $L^1(X)$. Then $P_{X}''\in\B(M(K))$ can be identified with $P_{\widetilde{X}}\in \B(M(K))$ for some measurable set $\widetilde{X}\subset K$. The collection
$\{\widetilde{X}: X \in \B_{\Omega}\}$ forms a base of clopen sets for the topology on $K$. Hence by Theorem \ref{4.6.17} we have the following. 

\begin{proposition}
Let $\Omega$ be a measure space, and take $1\leq q< \infty$. Then
$$
\norm{(\Phi_1,\ldots,\Phi_n)}^{(1,q)}_n=\sup_{\mathbf{X}}\lp \sum_{i=1}^n\norm{P_{X_i}''(\Phi_i)}^q\rp^{1/q}\quad (\Phi_1,\ldots,\Phi_n \in L^{\, 1}(\Omega)'')\,,
$$
where the supremum is taken over all measurable partitions $\mathbf{X}=(X_1,\ldots,X_n)$ of $\Omega$.\qed
\end{proposition}

\subsection{Extensions of multi-norms}
Let $F$ be a Banach space, and let $E$ be a multi-normed space. For each $n \in \N$ we define a norm $\norm{\,\cdot\,}_n^{\B}$ on the space $F^n$ by setting
$$
\norm{(x_1,\ldots,x_n)}^{\B}_n=\sup_{U}\norm{(U(x_1),\ldots,U(x_n))}_n\quad (x_1,\ldots, x_n \in F)\,,
$$
where the supremum is taken over all $U\in \B(F, E)_{[1]}$. It is immediately checked that this defines a multi-norm over $F$, and that \begin{equation}\label{eqI}
\mathcal{M}(F, E)=\B(F, E)\quad {\rm with}\quad \norm{T}_{mb}=\norm{T}\quad (T \in \mathcal{M}(F, E))\,.
\end{equation}
Let $(\tnorm{\,\cdot\,}_n: n \in \N)$ be a multi-norm over $F$ such that \eqref{eqI} holds. Then it is clear that $\norm{\,\cdot\,}^{\B}_n\leq \tnorm{\,\cdot\,}_n\,\; (n \in \N)$.

\begin{definition}
Let $F$ be a Banach space, let $E$ be a multi-normed space. Then the multi-norm $(\norm{\,\cdot\,}_n^{\B}: n \in \N)$ described above is the \italic{extension} to $F$ of the multi-norm on $E$.
\end{definition}

\begin{example}[{\cite[Example 4.2]{DP08}}]\label{4.6.22}
Let $\Omega$ be a measure space, and take $1< p\leq q < \infty$. For each measurable partition $\mathbf{X}=(X_1,\ldots,X_n)$ of $\Omega$ we define
\begin{equation}
\norm{(f_1,\ldots,f_n)}_n^{[p,q]}=\sup_{\mathbf{X}}\left(\sum_{i=1}^{n}\norm{\chi_{X_i}f_i}_{p}^{q}\right)^{1/q}\quad (f_1,\ldots,f_n \in L^{p}(\Omega))\,.
\end{equation}
Then $\norm{\,\cdot\,}_n^{[p,q]}$ is a norm on $L^{p}(\Omega)^n$ and the family $(\norm{\,\cdot\,}_{n}^{[p,q]}: n \in \N)$ is a multi-norm over $L^{p}(\Omega)$ called the \italic{standard $(p,q)$-multi-norm}. In the same way as Proposition \ref{4.6.15} we can show that
$$
\norm{(f_1,\ldots, f_n)}_n^{[p,p]}=\norm{\mod{f_1}\vee\cdots \vee \mod{f_n}}_{p}\,.
$$
\end{example}

\subsubsection{The extension of the standard $(p,q)$-multi-norm}
Let $F$ be a Banach space, and let $1< p\leq q < \infty$. Let $\Omega$ be an infinite locally compact space, and let $m$ be a positive, regular Borel measure on $\Omega$. We shall identify the extension to $F$ of the standard-$(p,q)$ multi-norm on $L^{p}(\Omega)=L^{p}(\Omega, m)$. We shall denote this multi-norm by $(\norm{\,\cdot\,}_n^{\B}: n \in \N)$.

Let $D_n$ denote the collection of all $(f_1,\ldots,f_n) \in (L^{p'}(\Omega)_{[1]})^n$ such that the sets $\supp f_1,\ldots, \supp f_n$ are pairwise disjoint. For a Banach space $X$, set
$$
B_n(X)=\lb (U(f_1),\ldots, U(f_n)): U \in \B(L^{p'}(\Omega), X)_{[1]},\, (f_1,\ldots, f_n) \in  D_n\rb\subset X^n\,.
$$

\begin{lemma}\label{4.6.19}
We have $B_n(X)=\lb x \in X^n: \mu_{{p}, n}(x)\leq 1\rb$.
\end{lemma}
\begin{proof}
Set $C_n(X)=\lb x \in X^n: \mu_{p, n}(x)\leq 1\rb$.

Take $U \in \B(L^{p'}(\Omega), X)_{[1]}$ and ($f_1,\ldots, f_n) \in  D_n$. Set $X_i=\supp f_i\,\; (i \in \N_n)$ and set $x=(U(f_1),\ldots, U(f_n)) \in X^{n}$. For each $\lambda \in X'_{[1]}$ we have
\begin{equation*}
\begin{split}
\lp\sum_{i=1}^n\mod{\la U(f_i),\lambda\ra}^{p}\rp^{1/p}&=\lp\sum_{i=1}^n\mod{\la f_i, U'(\lambda)\ra}^{p}\rp^{1/p}\leq \lp\sum_{i=1}^n\norm{\chi_{X_i}U'(\lambda)}_{p}^{p}\rp^{1/p} \\
&= \norm{U'(\lambda)}_{p}\leq 1\,.
\end{split}
\end{equation*}
Hence $\mu_{p,n}(\lambda)\leq 1$, and so $B_n(X)\subset C_n(X)$.

Conversely, take $x=(x_1,\ldots,x_n) \in C_n(X)$. Choose non-null, pairwise disjoint subsets $X_1,\ldots, X_n\subset \Omega$ with $m(X_i)<\infty\,\; (i \in \N_n)$ (this is possible because of our hypotheses on $\Omega$ and $m$). Set $f_i=\frac{\chi_{X_i}}{m(X_i)^{1/p'}}\,\; (i \in \N_n)$, so that $(f_1,\ldots,f_n) \in D_n$.

Set
$$
U=\sum_{i=1}^{n}x_i\otimes\frac{\chi_{X_i}}{m(X_i)^{1/p}}\in \B(L^{p'}(\Omega), X)\,.
$$
For $f \in L^{p'}(\Omega)$, we have
\begin{equation*}
\begin{split}
\norm{U(f)}&=\norm{\sum_{i=1}^n \la f, \frac{\chi_{X_i}}{m(X_i)^{1/p}}\ra x_i}\leq \mu_{p,n}(x)\lp \sum_{i=1}^n\mod{\la f, \frac{\chi_{X_i}}{m(X_i)^{1/p}}\ra}^{p'}\rp^{1/p'} \\
&\leq\mu_{p,n}(x)\lp \sum_{i=1}^n\norm{\chi_{X_i}f}^{p'}\rp^{1/p'}\leq\mu_{p,n}(x)\norm{f} \\
\end{split}
\end{equation*}
It follows that $U \in \B(L^{p'}(\Omega), X)_{[1]}$. Since $x=(U(f_1),\ldots, U(f_n))$, we have $B_n(X)=C_n(X)$, as required.
\end{proof}

Let $F$ be a Banach space. Since $L^{p'}(\Omega)$ is reflexive, every $T \in \B(L^{p'}(\Omega), F')$ can be written as $T=U'$ where $U \in \B(F, L^{p}(\Omega))$. Hence
\begin{equation}\label{dualmap}
B_n(F')=\lb (U'(f_1),\ldots, U'(f_n)): U \in \B(F, L^{p}(\Omega))_{[1]},\, (f_1,\ldots, f_n) \in  D_n\rb\,.
\end{equation}
The following corollary follows easily by taking duals of operators and applying Lemma \ref{4.6.19}.

\begin{corollary}\label{4.6.20}
Let $F$ be a Banach space, and let $1< p\leq q < \infty$. Let $\Omega$ be an infinite locally compact space, and let $m$ be regular Borel measure on $\Omega$. Then the extension to $F$ of the standard-$(p,q)$ multi-norm on $L^{p}(\Omega, m)$ is the weak-$(p,q)$ multi-norm on $F$.\qed
\end{corollary}

The following lemma is needed for the applications in the next section.

\begin{lemma}\label{4.6.21}
Let $F$ be a Banach space, and let $1< p\leq q < \infty$. Let $\Omega$ be an infinite locally compact space, and let $m$ be regular Borel measure on $\Omega$. Then we have
$$
\norm{(\Phi_1,\ldots,\Phi_n)}_n^{\mathcal{B}}=\sup \norm{(U''(\Phi_1),\ldots,U''(\Phi_n))}_n^{[p,q]}
$$
where $\Phi_1,\ldots,\Phi_n \in F''$ and the supremum is taken over all $U \in \B(F, L^{p}(\Omega, m))_{[1]}$.
\end{lemma}
\begin{proof}
Take $\Phi=(\Phi_1,\ldots,\Phi_n)\in (F'')^n$. By Corollary \ref{4.6.20} and Lemma \ref{4.6.6} we have
$$
\norm{\Phi}_n^{\mathcal{B}}=\sup\lb\lp\sum_{i=1}^n\mod{\la\lambda_i,\Phi_i\ra}^{q}\rp^{1/q}: \lambda \in (F')^n,\, \mu_{p,n}(\lambda)\leq 1\rb\,.
$$
By \eqref{dualmap} and Lemma \ref{4.6.19} this is equal to
\begin{equation*}
\begin{split}
\norm{\Phi}_n^{\mathcal{B}}&=\sup\lb\lp\sum_{i=1}^n\mod{\la U'(f_i),\Phi_i\ra}^{q}\rp^{1/q}: U \in \B(F, L^{p}(\Omega))_{[1]},\, (f_i) \in  D_n\rb \\
&=\sup\lb\lp\sum_{i=1}^n\mod{\la U''(\Phi_i),f_i\ra}^{q}\rp^{1/q}: U \in \B(F, L^{p}(\Omega))_{[1]},\, (f_i) \in  D_n\rb \\
&=\sup \lb\norm{(U''(\Phi_1),\ldots,U''(\Phi_n))}_n^{[p,q]}: U \in \B(F, L^{p}(\Omega))_{[1]}\rb\,.\qedhere
\end{split}
\end{equation*}
\end{proof}

\section{Generalized notions of amenability}
Let $G$ be a locally compact group with left Haar measure $m$, and let $L^1(G)=L^1(G,m)$. For $f \in L^1(G)$ and $s \in G$ we define $s\cdot f \in L^1(G)$ by $(s\cdot f)(t)=f(s^{-1}t)\,\; (t \in G)$. This defines an action of $G$ on the space $L^1(G)$. We can extend this action by duality to $L^1(G)''$. An element $\Lambda \in L^1(G)''$ is a {\it mean} if $1=\la1, \Lambda\ra=\norm{\Lambda}=1$, and {\it left invariant} if $\{s\cdot \Lambda:s \in G\}=\{\Lambda\}$. If there exists a left invariant mean $\Lambda \in L^1(G)''$, then $G$ is {\it amenable}.

Now we show how to use multi-norms to generalize this concept. The idea of using multi-norms in this way is due to Dales and Polyakov.

\begin{definition}
Let $G$ be a locally compact group, and take $1\leq p\leq q<\infty$. A mean $\Lambda \in L^{\, 1}(G)''$ is {\it left} $(p,q)$-{\it multi-invariant} if the set $\{s\cdot \Lambda : s \in G\}$ is multi-bounded in the weak $(p,q)$-multi-norm. If there exists such an element $\Lambda \in L^{\, 1}(G)''$, then $G$ is $(p,q)$-{\it amenable}. For such an element $\Lambda$, we set
$$
C_{p,q}(\Lambda)=\sup\{\norm{(s_1\cdot\Lambda, \ldots, s_n\cdot\Lambda)}^{(p,q)}_{n}: s_{1},\ldots,s_{n} \in G,\, n \in \N\}\,.
$$
\end{definition}

The idea behind this definition is an attempt to measure the `left-invariance' of a mean $\Lambda \in L^1(G)''$ by measuring the growth of the sets $\{s\cdot \Lambda: s \in F\}$ as $F$ ranges through all finite subsets of $G$. The following implications follow immediately from Corollary \ref{4.6.15a}:
$$(q,q)\text{-amenable} \Rightarrow (p,q)\text{-amenable} \Rightarrow (1,q)\text{-amenable}\quad (1< p < q<\infty)\,;$$
$$(p,p)\text{-amenable} \Rightarrow (q,q)\text{-amenable for all } q\geq p\,.$$
The strongest of these conditions is $(1,1)$-amenability. It follows from the multi-norm axiom (A4) that an amenable locally compact group is $(1,1)$-amenable. The converse will be shown in Proposition \ref{5.1.9}.

There is an obvious definition of a {\it right} $(p,q)$-{\it multi-invariant mean}. Set $A=L^{\, 1}(G)$. Let $\Lambda \in A''$ be a left $(p,q)$-multi-invariant mean. Define $T: A\rightarrow A$ by 
$$
T(a)(s)=a(s^{-1})\Delta(s^{-1})\quad (a \in A, s \in G)\,,
$$
where $\Delta$ is the modular function of $G$. Then $T'':A''\rightarrow A''$ takes the set $\{s\cdot \Lambda: s \in G\}$ to $\{T''(\Lambda)\cdot s: s \in G\}$, and $T'(1)=1$. By Proposition \ref{4.6.7}, $T'' \in \mathcal{M}(A'', A'')$, and hence $T''(\Lambda)$ is a right $(p,q)$-multi-invariant mean on $G$.

Of most interest to us are $(1,q)$-amenability and $(q,q)$-amenability, of which the latter concept is formally stronger.

\subsection{$(1,q)$-amenability}
Since any $(p,q)$-amenable group is $(1,q)$-amenable, it is particularly interesting to investigate this concept.

\begin{proposition}
Let $G$ be a locally compact group, and take $1\leq q<\infty$. Then there exists a $(1,q)$-multi-invariant mean in $L^{\, 1}(G)$ if and only if $G$ is compact.
\end{proposition}
\begin{proof}
If $G$ is compact, then $\Lambda=\chi_{G}/m(G)$ is an invariant and hence $(1,p)$-multi-invariant mean.

Conversely, assume towards a contradiction that $G$ is not compact and that there exists a $(1,q)$-multi-invariant mean $a \in L^{\, 1}(G)$. There is a compact set $V$ such that $c=\int_{V}\mod{a(t)} \dd m(t)\neq 0$. Since $G$ is not compact, for each $N \in \N$, there exist elements $s_1,\ldots,s_N \in G$ such that the sets $s_1V,\ldots,s_NV$ are pairwise disjoint. We have $\chi_{s_{i}V}(s_{i}\cdot a)=s_{i}\cdot(\chi_{V}a)$ and so
$$
C_{1,q}(a)\geq \norm{(s_1\cdot a,\cdots,s_n\cdot a)}_n^{(1,q)}\geq\lp\sum_{i=1}^{N}\norm{\chi_{s_{i}V}(s_{i}\cdot a)}^{q}\rp^{1/q}=\lp\sum_{i=1}^{N}\norm{\chi_{V}a}^{q}\rp^{1/q}=N^{1/q}c\,.
$$
This holds for all $N \in \N$, the required contradiction.
\end{proof}

The following result was first proved by Dales and Polyakov for discrete groups.

\begin{proposition}\label{5.1.9} \ifmargin \marginpar[5.1.9]{5.1.9} \fi
Let $G$ be a locally compact group. Then $G$ is $(1,1)$-amenable if and only if $G$ is amenable.
\end{proposition}
\begin{proof}
It is clear that every amenable locally compact group is $(1,1)$-amenable.

We set $A=L^{\, 1}(G)$. Suppose that $G$ is $(1,1)$-amenable, and let $\Lambda \in L^1(G)''$ be a $(1,1)$-multi-invariant mean. For each $U \in \B_G$ we define $\la \chi_{U}, \widetilde{\Lambda}\ra \in \R^{+}$ by
$$
\la \chi_{U}, \widetilde{\Lambda}\ra=\sup \sum_{i=1}^n\la \chi_{X_i}, s_i\cdot \Lambda\ra\,,
$$
where the supremum is taken over all $n \in \N$, all $s_1,\ldots, s_n \in G$, and all measurable partitions $(X_1,\ldots,X_n)$ of $U$. The supremum is finite since $\Lambda$ is $(1,1)$-invariant. Since
$$
\la \chi_{U\cup V}, \widetilde{\Lambda}\ra=\la \chi_{U}, \widetilde{\Lambda}\ra+\la \chi_{V}, \widetilde{\Lambda}\ra\,,
$$
for all $U, V \in \B_{G}$ with $U\cap V=\emptyset$, we can extend $\widetilde{\Lambda}$ to a linear map $\widetilde{\Lambda}:\mathcal{S}=\lin \lb \chi_{U}: U \in \B_G\rb\rightarrow \C$ by setting
$$
\la \sum_{i=1}^n\alpha_i\chi_{U_i},\widetilde{\Lambda}\ra=\sum_{i=1}^n\alpha_i\la\chi_{U_i},\widetilde{\Lambda}\ra\,.
$$
The set $\mathcal{S}$ is dense in $A'$, and $\mod{\la \lambda, \widetilde{\Lambda}\ra}\leq C_{1,1}(\Lambda)\norm{\lambda}_{\infty}\,\; (\lambda \in \mathcal{S})$. Hence $\widetilde{\Lambda}$ extends to an element $\widetilde{\Lambda} \in A''$ with $\norm{\widetilde{\Lambda}}\leq C_{1,1}(\Lambda)$. It is easily checked that $s\cdot\widetilde{\Lambda}=\widetilde{\Lambda}\,\; (s \in G)$, and $\la 1, \widetilde{\Lambda}\ra\geq \la 1,\Lambda\ra=1$. This implies that $G$ is amenable.
\end{proof}

\begin{proposition}
The free group on two generators is not $(1,q)$-amenable for any $q\geq 1$.
\end{proposition}
\begin{proof}
Let $\F_2$ denote the free group on the generators $a, b$. Then each element of $\F_2$ is a reduced word in the alphabet $\{a, a^{-1}, b, b^{-1}\}$. For each $x \in \{a, a^{-1}, b, b^{-1}\}$ we set
$$
W(x)=\{w \in \F_2: w \text{ starts with } x\}\,,
$$
so that $\F_2$ is a disjoint union $\F_2=\{e\}\cup W(a)\cup W(a^{-1})\cup W(b)\cup W(b^{-1})$.

Assume towards a contradiction that $\F_2$ is $(1,q)$-amenable, and let $\Lambda \in L^1(G)''$ be a $(1,q)$-multi-invariant mean. Since
$$
1=\la 1, \Lambda\ra=\la \delta_{e}, \Lambda\ra+\la \chi_{W(a)}, \Lambda\ra+\la \chi_{W(a^{-1})}, \Lambda\ra+\la \chi_{W(b)}, \Lambda\ra+\la \chi_{W(b^{-1})}, \Lambda\ra\,,
$$
we have $0<\la\chi_{X},\Lambda \ra$ for some set $X\in \{\{e\}, W(a), W(a^{-1}), W(b), W(b^{-1})\}$. We may suppose that $\la \chi_{W(a)}, \Lambda\ra>0$. For each $n \in \N$ the sets $bW(a),\ldots,b^nW(a)$ are pairwise disjoint. Hence we have
$$
C_{1,q}(\Lambda)\geq \norm{(b\cdot \Lambda,\ldots, b^n\cdot\Lambda)}_n^{(1,q)}\geq \lp\sum_{i=1}^n\mod{\la\chi_{b^iW(a)},b^i\cdot\Lambda \ra}^q \rp^{1/q}=\la \chi_{W(a)}, \Lambda\ra n^{1/q}\,.
$$
This is true for each $n \in \N$, which is a contradiction.
\end{proof}

\subsection{Pseudo-amenability}
Here we show that $(1,q)$-amenability implies pseudo-amenability.

For a locally compact group $G$, we set
$$
P(G)=\lb f \in L^1(G): f\geq 0,\, \norm{f}=1\rb\,.
$$

\begin{proposition}\label{5.1.11} \ifmargin \marginpar[5.1.11]{5.1.11} \fi
Let $G$ be a locally compact group, and take $1\leq p\leq q<\infty$. Suppose that $G$ is $(p,q)$-amenable. Then there exists $C\geq 1$ such that, for each $n \in \N$, and for each finite set $\{s_{1},\ldots, s_{n}\} \subset G$, there exists $a \in P(G)$ with
$$
\norm{(s_{1}\cdot a,\ldots,s_{n}\cdot a)}_{n}^{(p,q)}\leq C\,.
$$
\end{proposition}
\begin{proof}
We set $A=L^{\, 1}(G)$.

Let $\Lambda \in A''$ be a $(p,q)$-multi-invariant mean. Set $C=C_{p,q}(\Lambda)+1$. Fix $n \in \N$ and a finite set $\{s_{1},\ldots, s_{n}\} \subset G$. By \cite[Proposition (0.1)]{PA} there is a net $(a_{\alpha})$ in $P(G)$ such that $\lim_{\alpha} a_{\alpha}=\Lambda$ in the weak-$*$ topology on $A''$. By Theorem \ref{duality} $(A^n, \norm{\,\cdot\,}^{(p,q)}_n)''=((A'')^n,\norm{\,\cdot\,}^{(p,q)}_n)$. Hence there is a net $b_{\alpha}=(b_{1, \alpha},\ldots,b_{n, \alpha})$ in $A^n$ such that
$$
\lim_{\alpha} b_{\alpha}=(s_{1}\cdot \Lambda,\ldots, s_{n}\cdot\Lambda)
$$
in the weak-$*$ topology on $(A'')^{n}$ and
$$
\sup_{\alpha}\norm{b_{\alpha}}_n^{(p,q)}\leq \norm{(s_{1}\cdot \Lambda,\ldots, s_{n}\cdot\Lambda)}_n^{(p,q)}\leq c_{\Lambda}\,.
$$
We can suppose that these nets are indexed by the same directed set. We have
$$
\lim_{\alpha} (s_{1}\cdot a_{\alpha}-b_{1, \alpha},\ldots, s_{1}\cdot a_{\alpha}-b_{n, \alpha})=0
$$
in the weak topology on $A^{n}$. By Mazur's theorem there is some convex combination
\begin{equation*}
\begin{split}
v&=\sum_{j=1}^{N}t_{j}(s_{1}\cdot a_{\alpha_{j}}-b_{1, \alpha_{j}},\ldots,s_{n}\cdot a_{\alpha_{j}}-b_{n, \alpha_{j}})\\
&=\lp s_{1}\cdot \lp \sum_{j=1}^{N}t_{j}a_{\alpha_{j}}\rp-\sum_{j=1}^{N}t_{j}b_{1,\alpha_{j}},\ldots,s_{n}\cdot\lp\sum_{j=1}^{N}t_{j} a_{\alpha_{j}}\rp-\sum_{j=1}^{N}t_{j}b_{n, \alpha_{j}}\rp
\end{split}
\end{equation*}
such that $\norm{v}^{(p,q)}_n\leq 1$. Set $a=\sum_{j}t_{j}a_{\alpha_j}$ and $b_i=\sum_{j}t_{j}b_{i, \alpha_j}\,\; (i \in \N_n)$.

Then we have
\begin{equation*}
\begin{split}
\norm{(s_1\cdot a,\ldots,s_n\cdot a)}_n^{(p,q)} &\leq \norm{v}_n^{(p,q)}+\norm{(b_1,\ldots,b_n)}_n^{(p,q)} \\
& \leq 1+C_{p,q}(\Lambda)= C\,.\qedhere
\end{split}
\end{equation*}
\end{proof}

Let $S$ be a set, and let $n \in \N$. Then $\mathcal{P}_{n}(S)$ denotes the collection of subsets of $S$ containing $n$ elements.

\begin{lemma}\label{5.1.12} \ifmargin \marginpar[5.1.12]{5.1.12} \fi
Let $f=\sum_{k=1}^N\beta_k\chi_{S_k} \in P(G)$ where $S_1\subset S_2 \subset \cdots \subset S_N \subset G$. Let $F=\{s_1,\ldots, s_n\} \in \mathcal{P}_{n}(G)$. Then
$$
\norm{(s_1\cdot f,\ldots,s_n\cdot f)}_n^{(1,1)}=\sum_{k=1}^N\mod{\beta_k}m(FS_k)\,.
$$
\end{lemma}
\begin{proof}
For $s \in G$ and $i \in \N_n$ we define $k(s, i) \in \N_N$ by
$$
k(s,i)=\min \lb k \in \N_N: s \in s_iS_k\rb\,.
$$
Also for $s \in G$ we define $k(s) \in \N_N$ by
$$
k(s)=\min \lb k(s,i): i \in \N_n\rb\,.
$$
Now we have
$$
(s_i\cdot f)(t)=\sum_{k=1}^n\beta_{k}\chi_{s_iS_k}(t)=\sum_{k=k(t,i)}^N\beta_k\quad (i \in \N_n,\, t \in G)\,,
$$
and
$$
\max_{i \in \N_n} \mod{s_i\cdot f}(t)=\sum_{k=k(t)}^N\mod{\beta_k}=\sum_{k=1}^N\mod{\beta_k}\chi_{\{s:k(s)\leq k\}}(t)\quad (t \in G)\,.
$$
Now, for each $k \in \N_N$, we have
\begin{equation*}
\begin{split}
\{s:k(s)\leq k\}&=\{s: \exists i \in \N_n,\, k(s, i)\leq k\} \\
&=\{s: \exists i \in \N_n, \exists l\leq k,\, s \in s_iS_l\} \\
&=FS_k\,.
\end{split}
\end{equation*}
Hence
\begin{equation*}
\norm{(s_1\cdot f,\ldots,s_n\cdot f)}_n^{(1,1)}=\int_{G}\max_{i \in \N_n} \mod{s_i\cdot f}(t) \dd m(t)=\sum_{k=1}^N\mod{\beta_k}m(FS_k)\,.\qedhere
\end{equation*}
\end{proof}

For a discrete group $G$ the condition in the next proposition is the same as the condition arrived at in \cite[Proposition 5.11]{DP}.

\begin{proposition}\label{5.1.13} \ifmargin \marginpar[5.1.13]{5.1.13} \fi
Let $G$ be a locally compact group, and take $1<q<\infty$. Suppose that $G$ is $(1,q)$-amenable. Then there exists $C\geq 1$ such that, for every $n \in \N$ and every $F \in \mathcal{P}_{n}(G)$, there exists a non-null, compact subset $S \subset G$ with
$$
\frac{m(FS)}{m(S)}\leq Cn^{1-1/q}\,.
$$
\end{proposition}
\begin{proof}
We set $A=L^{\, 1}(G)$. Let $C_0$ be the constant given in Proposition \ref{5.1.11}. Take $n \in \N$ and $F=\{s_1,\ldots,s_n\} \in \mathcal{P}_n(G)$. By Proposition \ref{5.1.11} there exists $a \in P(G)$ with $$\norm{(s_1\cdot a, \ldots,s_n\cdot a)}_n^{(1,1)}\leq n^{1/q'}\norm{(s_1\cdot a, \ldots,s_n\cdot a)}_n^{(1,q)}\leq C_0n^{1/q'}\,.$$

There exists $f \in P(G)$ with finite range, such that $\norm{f-a}\leq n^{1/q'-1}$. Then we have
\begin{equation*}
\begin{split}
\norm{(s_1\cdot f, \ldots,s_n\cdot f)}_n^{(1,1)}&\leq \norm{(s_1\cdot a, \ldots,s_n\cdot a)}_n^{(1,1)}+\norm{(s_1\cdot (f-a), \ldots,s_n\cdot (f-a))}_n^{(1,1)} \\
&\leq C_0n^{1/q'}+\sum_{i=1}^n\norm{s_i\cdot(f-a)}\leq (C_0+1)n^{1/q'}\,.
\end{split}
\end{equation*}

Set $C=C_0+1$. We can write 
$$
f = \sum_{k=1}^{N} \alpha_k \chi_{S_k}/m(S_k)\,,
$$
where $S_1\subset S_2 \subset \cdots \subset S_N \subset G$, where $\alpha_1,\dots, \alpha_N > 0$, and where $\sum_{k=1}^{N} \alpha_k =1$. By Lemma \ref{5.1.12} we have
$$
\sum_{k=1}^N\frac{\alpha_k m(FS_k)}{m(S_k)}\leq Cn^{1/q'}\,.
$$
The left-hand side of this inequality is a convex sum, hence there exists $k\in \N_N$ such that
$$
\frac{m(FS_k)}{m(S_k)}\leq Cn^{1/q'}\,.
$$
Finally we set $S=S_k$, giving the result.
\end{proof}

A discrete group $G$ satisfying the condition in the next proposition is called \italic{pseudo-amenable} in \cite[Definition 5.5]{DP}.

\begin{proposition}
Let $G$ be a locally compact group, and take $1<q<\infty$. Suppose that $G$ is $(1,q)$-amenable. Then for all $\varepsilon >0$, there exists $n_0=n_0(\varepsilon) \in \N$ such that, for all $n \geq n_0$ and $F \in \mathcal{P}_{n}(G)$, there exists a non-null, compact subset $S \subset G$ with
$$
\frac{m(FS)}{m(S)}\leq \varepsilon n\,.
$$
\end{proposition}
\begin{proof}
Let $C\geq 1$ be the constant prescribed in Proposition \ref{5.1.13}. Take $\varepsilon >0$, and choose $n_0 \in \N$
with $n_0^{1/q}=n_0^{1-1/q'}\geq C/\varepsilon$, so that $Cn^{1/q'}\leq \varepsilon n\,\; (n\geq n_0)$. Now take $n \in \N$ and $F \in \mathcal{P}_n(G)$. By Proposition \ref{5.1.13} there exists a non-null, compact subset $S \subset G$ with $\frac{m(FS)}{m(S)}\leq Cn^{1/q'}$. Hence $\frac{m(FS)}{m(S)}\leq n\varepsilon$, as required.
\end{proof}

\begin{remark}
The following facts are proved in \cite{DP}:\smallskip

{\rm (i)} Every subgroup of a pseudo-amenable discrete group is pseudo-amenable.\smallskip

{\rm (ii)} The free group on $2$ generators, $\F_2$ is not pseudo-amenable.
\end{remark}

\section{Injectivity of the $L^1(G)$-module $L^p(G)$}
Let $G$ be a locally compact group. We now consider $L^1(G)$ as a Banach algebra equipped with the {\it convolution product} $\star$ given by
\begin{equation}\label{convolution}
(f \star g)(s)=\int_{G} f(t)g(t^{-1}s) \dd m(t)\quad (s \in G)\,,
\end{equation}
where $f, g \in L^1(G)$, and the integral is defined for almost all $s \in G$.

We denote by $\varphi_G$ the {\it augmentation character} on $G$, given by
$$
\varphi_G(f)=\int_G f(t) \dd m(t)\quad (f \in L^1(G))\,.
$$

Take $1< p<\infty$, and let $L^p(G)=L^{\, p}(G, m)$. Let $f \in L^1(G)$, and let $g \in L^p(G)$. Then again we can define $f \star g$ on $G$ via \eqref{convolution} and we have $f\star g \in L^p(G)$. With this multiplication $L^p(G)$ has the structure of a Banach left $L^1(G)$-module, and further $L^p(G)$ is a faithful $L^1(G)$-module (see \cite[Theorem 3.3.19]{HGD}).

It is convenient to write the module multiplication as a Banach space valued integral. For each $f \in L^{1}(G)$ and $g \in L^p(G)$ we have
\begin{equation}\label{BSint}
f\cdot g=\int f(t)t\cdot g \dd m(t)\,.
\end{equation}
This is a special case of \cite[Proposition 2.1]{JO72}.

\subsection{A coretraction problem}
Again, let $G$ be a locally compact group, and take $1<p<\infty$. We set $A=L^{1}(G)$, $E=L^{p}(G)$ and $J=\mathcal{B}(A,E)$. We now define an action of $G$ on the space $J$ by 
\begin{equation}
\left(t * U\right)(a)=t\cdot U(t^{-1} \cdot a)\quad (a \in A)\,.
\end{equation}
For each $U \in J$ and $a \in A$, the map $t\mapsto (t*U)(a)=t\cdot U(t^{-1}\cdot a),\, G\rightarrow E$ is continuous. This follows from the inequality
\begin{equation*}
\begin{split}
\left\|t\cdot U(t^{-1}\cdot a)-U(a)\right\|_{p}&\leq\left\|t\cdot U(t^{-1}\cdot a)-t\cdot U(a)\right\|_{p}+\left\|t\cdot U(a)-U(a)\right\|_{p}\\
                         &=\left\|U(t^{-1}\cdot a-a)\right\|_{p}+\left\|t\cdot U(a)-U(a)\right\|_{p} \\
                         &\leq \norm{U}\norm{t^{-1}\cdot a-a}_{1}+\left\|t\cdot U(a)-U(a)\right\|_{p}\,,\end{split}
\end{equation*}
and \cite[3.3.11]{HGD}.
\begin{proposition}
There is a Banach left $A$-module structure on $J$ given by
\begin{equation}
\begin{split}
\left(b * U\right)(a)&=\int_{G}b(t)\left(t * U\right)(a) \dd m(t) \quad (a, b \in A,\, U \in J)\,.\end{split}
\end{equation}
\end{proposition}
\begin{proof}
This is similar to the standard proof that $L^{p}(G)$ is a left $L^{1}(G)$-module \cite[3.3.19]{HGD}. Fix $U \in J$ and $a, b \in A$. Let $\psi \in C_{00}(G)$. By \holders inequality, we have
$$
\int_{G}\left|U(t^{-1}\cdot a)(t^{-1}s)\right|\left|\psi(s)\right| \dd m(s)\leq \norm{t\cdot U(t^{-1}\cdot a)}_{p}\norm{\psi}_{p'}\leq \norm{U}\norm{a}_{1}\norm{\psi}_{p'}
$$
for each $t \in G$. Now define
\begin{equation*}
\begin{split}
\Lambda: C_{00}(G)\rightarrow \mathbb{C}: \psi\mapsto &\int_{G}(b * U)(a)(s)\psi(s) \dd m(s) \\
                         &=\int_{G}\left(\int_{G}b(t)U(t^{-1}\cdot a)(t^{-1}s) \dd m(t)\right)\psi(s) \dd m(s)\\
                         &=\int_{G}b(t)\left(\int_{G}U(t^{-1}\cdot a)(t^{-1}s)\psi(s) \dd m(s)\right) \dd m(t)\,.
\end{split}
\end{equation*}
Then $\left|\Lambda(\psi)\right|\leq \norm{b}_{1}\norm{U}\norm{a}_{1}\norm{\psi}_{p'}\,\; (\psi \in C_{00})$, and so $\Lambda$ extends to an element of $L^{p'}(G)'$ of norm at most $\norm{b}_{1}\norm{U}\norm{a}_{1}$. Hence by the identification of $L^{p'}(G)'$ with $L^{p}(G)$, we have $(b*U)(a) \in L^{p}(G)$ and $b*U \in J$ with $\norm{b* U}\leq \norm{b}_{1}\norm{U}$.

Associativity of $*$ follows in the same way as \cite[Proposition 2.1]{JO72}.
\end{proof}

We shall denote this left $A$-module by $\widetilde{J}=(J, \,*\,)$. (We could similarly define a right multiplication such that $J$ becomes an $L^{1}(G)$-bimodule).

Now we define an embedding $\widetilde{\Pi} :E \rightarrow \widetilde{J}$, by
$$
(\widetilde{\Pi}x)(a)=\varphi_{G}(a)x\quad (a \in A)\,.
$$
For $b \in A$, we have
$$
\left(b * \widetilde{\Pi}x\right)(a)=\int_G b(t)\varphi_G(t^{-1}\cdot a)t\cdot x \dd m(t)=\varphi_G(a)b\star x=\widetilde{\Pi}(b\star x)(a)\quad (a \in A)\,,
$$
and so $\widetilde{\Pi}$ is a left $A$-module morphism; further, $\widetilde{\Pi}$ is admissible (a splitting operator is $U \mapsto U(a_{0})$ for any $a_{0} \in A$ with $\varphi_{G}(a_{0})=1$).

\begin{proposition}\label{coretract} \ifmargin \marginpar[coretract]{coretract} \fi
Let $G$ be a locally compact group, and let $1<p<\infty$. Then $L^{p}(G)$ is injective in $\Lmod$ if and only if the morphism $\widetilde{\Pi}$ is a coretraction in $\Lmod$.
\end{proposition}
\begin{proof}
It is clear that the condition is necessary, we shall prove sufficiency.

As above we set $A=L^1(G)$, $E=L^p(G)$, and $J=\B(A, E)$. Also set $F=L^{p'}(G)$.

Suppose that $\widetilde{\Pi}$ is a coretraction, so that there exists $R \in \AB(\widetilde{J}, E)$ with $R\circ \widetilde{\Pi}=I_{E}$. For $f \in A$ we define $Q_{f} \in \B(J)$ by
$$
Q_{f}(U)(a)=(a * U)(f)\quad (a \in A,\, U \in J)\,.
$$
For $x \in E$ and $a \in A$, we have
\begin{equation}\label{Q1}
Q_{f}\lp\Pi(x)\rp(a)=(a * \Pi(x))(f)=\varphi_{G}(a)(f \cdot x)=\lp\widetilde{\Pi} (f\cdot x)\rp(a)\,,
\end{equation}
and for $U \in J$, and $a,b \in A$ we have
\begin{equation*}
\begin{split}
(b * Q_{f}(U))(a)&=\int_{G}b(s)s \cdot Q_{f}(U)(s^{-1}\cdot a)\dd m(s)\\
												&=\int_{G}b(s)s \cdot((s^{-1}\cdot a)*U)(f)\dd m(s) \\
                        &=\int_{G}b(s)s \cdot\left(\int_{G}a(st)t\cdot U(t^{-1}\cdot f)\dd m(t)\right)\dd m(s)\\
                        &=\int_{G}b(s)s \cdot\left(\int_{G}a(t)s^{-1}t\cdot U(t^{-1}s\cdot f)\dd m(t)\right)\dd m(s)\\
                        &=\int_{G}b(s)\left(\int_{G}a(t)t\cdot U(t^{-1}s\cdot f)\dd m(t)\right)\dd m(s)\quad \text{(by \cite[III.6.20]{DS})}\\
                        &=\int_{G}a(t)\left(\int_{G}b(s)t\cdot U(t^{-1}s\cdot f)\dd m(s)\right) \dd m(t)\quad \text{(by Fubini)} \\
                        &=\int_{G}a(t)t\cdot U\left(\int_{G}b(s)t^{-1}s\cdot f\dd m(s)\right) \dd m(t) \quad \text{(by \cite[III.6.20]{DS})} \\
                        &=\int_{G}a(t)t\cdot U(t^{-1}\cdot b\star f) \dd m(t) \quad \text{(by \eqref{BSint})}\\
                        &=(a *U)(b\star f)\,.\end{split}
\end{equation*}
Hence
\begin{equation}\label{Q2}
(b * Q_{f}(U))(a)=(a *U)(b\star f)\,.
\end{equation}
We also have
\begin{equation}\label{Q3}
Q_f(b\cdot U)(a)=\int_G a(t)t\cdot U(t^{-1}\cdot f\star b)=(a *U)(f\star b)\,.
\end{equation}

Let $(e_{\alpha})$ be a bounded approximate identity for $A$, and set $Q_{\alpha}=Q_{e_{\alpha}}$. Let $Q$ be a weak-$*$ cluster point in $\B(J, J)=\lp J\tensor (A\tensor F)\rp'$ of the bounded net $(Q_{\alpha})$. By passing to a subnet we may suppose that $Q=\lim_{\alpha} Q_{\alpha}$ in the weak-$*$ topology. Take $x \in E$. Then for each $a \in A$ and $\lambda \in F$, by \eqref{Q1} we have
$$
\la \lambda, Q(\Pi x)(a)\ra=\lim_{\alpha}\la \lambda, \lp\widetilde{\Pi} (e_{\alpha}\cdot x)\rp(a)\ra=\la \lambda, \lp\widetilde{\Pi} x\rp(a)\ra\,.
$$
Hence $Q\circ \Pi=\widetilde{\Pi}$. Take $U \in J$ and $b \in A$. Then for each $a \in A$, and $\lambda \in F$, by \eqref{Q2} and \eqref{Q3} we have
\begin{equation*}
\begin{split}
\la \lambda, \lp b*Q(U)\rp(a)\ra &=\lim_{\alpha}\la \lambda, (a *U)(b\star e_{\alpha})\ra=\lim_{\alpha}\la \lambda, (a *U)(e_{\alpha}\star b)\ra \\
&=\lim_{\alpha}\la \lambda, Q_{\alpha}(b\cdot U)(a)\ra=\la \lambda, Q(b\cdot U)\ra\,.
\end{split}
\end{equation*}
Hence $b*Q(U)=Q(b\cdot U)$ and $Q \in \AB(J, \widetilde{J})$.

Finally we set $\rho=R \circ Q$, then $\rho \in \AB(J,E)$ and $\rho\circ \Pi=I_{E}$. Therefore $E$ is injective in $\Amod$.
\end{proof}
\subsection{Main result}
Let $G$ be a locally compact group, and let $1<p<\infty$. We shall prove that, if $L^{p}(G)$ is injective in $\Lmod$, then $G$ must be $(p,p)$-amenable.

We start with a generalization of \cite[Lemma 5.2]{DP}. For $n \in \N$, we set $D_n=\{-1, 1\}^n$,
and for $j \in \N_n$ we set
$$
D^+_n(j)=\{(d_{1},\ldots, d_{n}) \in D_n: d_j=1\},\quad D^-_n(j)=\{(d_{1},\ldots, d_{n}) \in D_n: d_j=-1\}\,.
$$
\begin{lemma}\label{5.2.3} \ifmargin \marginpar[5.2.3]{5.2.3} \fi
Let $n\in \mathbb{N}$, let $E$ be a normed space, and let $F: \mathbb{N}_{n}\times \N_{n}\rightarrow E$. Set 
$$
C= {\rm max}\left\{\left(\sum_{j=1}^{n}\norm{\sum_{i=1}^{n}d_{i}F(i, j)}^{p}\right)^{1/p}:(d_1,\dots,d_n)\in D_n\right\}\,.
$$
Then
$$
\left(\sum_{j=1}^{n}\norm{F(j,j)}^p\right)^{1/p}\leq C\,.
$$
\end{lemma}
\begin{proof} 
Take $d=(d_1,\dots,d_n)\in D_n$, and set $x_{j,d}=\sum^n_{i=1}d_i F(i,j)\;\,(j\in \N_n)$. By hypothesis, we have $\sum^n_{j=1}\norm{x_{j,d}}^p\leq C^{p}$. Since there are $2^n$ elements in $D_n$, we have
$$
\sum^n_{j=1}\sum_{d\in D_n}\norm{x_{j,d}}^p \leq 2^nC^{p}\,.
$$
For each $j \in \N_n$ we can write the term $ \sum_{d\in D_n}\norm{x_{j,d}}^p$ as
\begin{equation*}
\begin{split}
\sum_{d\in D_n}\norm{x_{j,d}}^p &=\sum_{d \in D^+_{n}(j)}\norm{x_{j,d}}^{p}+\sum_{d \in D^-_{n}(j)}\norm{x_{j,d}}^{p}  \\
           &= \sum_{d \in D_{n-1}}\norm{\sum_{i\neq j}d_{i}F(i,j)+F(j,j)}^{p}+ \norm{\sum_{i\neq j}d_{i}F(i,j)-F(j,j)}^{p}\\
           &\geq \sum_{d \in D_{n-1}}2\norm{F(j,j)}^{p} \quad \text{(by Jensen's inequality)} \\
           &=2^{n-1}.2\norm{F(j,j)}^{p}=2^{n}\norm{F(j,j)}^{p}\,.
\end{split}
\end{equation*}
This holds for each $j \in \N_n$, and so summing over $j$ we get
$$
2^n\sum^n_{j=1}\norm{F(j,j)}^p\leq\sum^n_{j=1}\sum_{d\in D_n}\norm{x_{j,d}}^p\leq 2^nC^p\,.
$$
Hence we have $\sum^n_{j=1}\norm{F(j,j)}^p\leq C^p$, and the result follows.
\end{proof}

For a measurable subset $V \subset G$ and $U \in J$ we define $\chi_{V}U \in J$ by
$$
(\chi_{V}U)(a)(s)=\chi_{V}(s)U(a)(s)\quad (a \in A,\, s \in G)\,.
$$

\begin{proposition}\label{5.3.4} \ifmargin \marginpar[5.3.4]{5.3.4} \fi
Let $\Omega$ be a measure space, and take $p$ with $1 < p <\infty$. Let $R:\B(L^{\, 1}(\Omega), L^{\,p}(\Omega))\rightarrow L^{\, p}(\Omega)$ be a bounded linear operator, and let $(X_{i})_{i=1}^{n}$ and $(Y_{i})_{i=1}^{n}$ be measurable partitions of $\Omega$. Then
$$
\left(\sum_{i=1}^{n}\norm{\chi_{X_{i}}R(\chi_{Y_{i}}U)}^{p}_{p}\right)^{1/p} \leq \norm{R}\norm{U}\quad (U \in \B(L^{\, 1}(\Omega), L^{\,p}(\Omega)))\,.
$$
\end{proposition}
\begin{proof}
Define $F:\N_{n}\times \N_{n}\rightarrow L^{\, p}(\Omega)$ by
$$
F(i, j)=\chi_{X_{j}}R(\chi_{Y_{i}}U)\quad (i, j \in \N_{n})\,.
$$
For each $(d_{1},\ldots, d_n) \in D_n$ we have
\begin{equation*}
\begin{split}
\sum_{j=1}^{n}\norm{\sum_{i=1}^{n}d_{i}F(i,j)}^{p}_p&=\sum_{j=1}^{n}\norm{\sum_{i=1}^{n}d_{i}\chi_{X_{j}}R(\chi_{Y_{i}}U)}^{p}_p=\sum_{j=1}^{n}\norm{\chi_{X_{j}}R\left(\sum_{i=1}^{n}d_{i}\chi_{Y_{i}}U\right)}^{p}_p \\
&=\norm{R\left(\sum_{i=1}^{n}d_{i}\chi_{Y_{i}}U\right)}^{p}_p\leq \norm{R}^{p}\norm{\sum_{i=1}^{n}d_{i}\chi_{Y_{i}}U}^{p}= \norm{R}^{p}\norm{U}^{p}\,.
\end{split}
\end{equation*}
Hence by Lemma \ref{5.2.3} we have
\begin{equation*}
\left(\sum_{j=1}^{n}\norm{F(j, j)}^{p}_p\right)^{1/p}=\left(\sum_{j=1}^{n}\norm{\chi_{X_{j}}R(\chi_{Y_{j}}U)}^{p}_p\right)^{1/p}\leq \norm{R}\norm{U}\,.\qedhere
\end{equation*}
\end{proof}

\begin{lemma}\label{5.3.5} \ifmargin \marginpar[5.3.5]{5.3.5} \fi
Let $U \in \B(L^{\, 1}(\Omega), L^{\, p}(\Omega))$, let $f_1,\ldots, f_n \in L^{\, p'}(\Omega)$ have pairwise disjoint supports, and let $x_1,\ldots, x_n \in L^{\, p}(\Omega)$ have pairwise disjoint supports. Set
$$
T=\sum_{i=1}^nx_i\otimes U'(f_i): L^{\, 1}(\Omega)\rightarrow L^{\, p}(\Omega)\,.
$$
Then $T \in \B(L^{\, 1}(\Omega), L^{\, p}(\Omega))$ and $\norm{T}\leq \norm{U}\max\{\norm{f_i}_{p'}\norm{x_i}_p: i \in \N_n\}$.
\end{lemma}
\begin{proof}
Set $X_i=\supp f_i\,\; (i \in \N_n)$, and set $C=\max\{\norm{f_i}_{p'}\norm{x_i}_p: i \in \N_n\}$. For $a \in L^{\, 1}(\Omega)$, we have
\begin{equation*}
\begin{split}
\norm{T(a)}_p^p&=\norm{\sum_{i=1}^n\la U(a), f_i\ra x_i}_p^p=\sum_{i=1}^n\mod{\la U(a), f_i\ra}^p \norm{x_i}_p^p \\
&\leq\sum_{i=1}^n\norm{\chi_{X_i}U(a)}^p_p\norm{f_i}_{p'}^p\norm{x_i}_p^p \\
&\leq C^p\norm{\chi_{X_1\cup\cdots\cup X_n}U(a)}_p^p\leq C^p\norm{U(a)}_p^p\,.
\end{split}
\end{equation*}
Therefore $\norm{T(a)}_p\leq C\norm{U(a)}_p$, and the result follows.
\end{proof}

\begin{lemma}\label{5.3.6} \ifmargin \marginpar[5.3.6]{5.3.6} \fi
Let $G$ be a locally compact group, and let $s_{1},\ldots, s_n \in G$. Then there exists an open, relatively compact neighbourhood $V$ of $e$ such that the sets $s_{1}V,\ldots,s_{n}V$ are pairwise disjoint.
\end{lemma}
\begin{proof}
Since $G$ is a Hausdorff space there exist pairwise disjoint open sets $U_{1},\ldots, U_{n}$ with $s_{i} \in U_{i}\,\; (i \in \N_n)$. For each $i \in \N_n$ the map $t\mapsto s_{i}t$ is continuous at $e$ and so there exists an open neighbourhood $V_{i}$ of $e$ with $s_{i}V_{i}\subset U_{i}\,\; (i \in \N_n)$. Set $V=\cap V_{i}$. Then $V$ is the required set.
\end{proof}

In the theorem below we shall use the following identity. For each $x \in L^{\, p}(G)$, $\lambda \in L^{\infty}(G)$ and $s \in G$ we have
\begin{equation}\label{eq4.3}
x\otimes (\lambda\cdot s)=s^{-1}*\left[(s\cdot x)\otimes \lambda\right]\,.
\end{equation}

\begin{theorem}\label{5.3.7} \ifmargin \marginpar[5.3.7]{5.3.7} \fi
Let $G$ be a locally compact group, and take $p$ with $1<p< \infty$. Suppose that $L^{\, p}(G)$ is injective in $\Lmod$. Then $G$ is $(p,p)$-amenable.
\end{theorem}
\begin{proof}
We set $A=L^{\, 1}(G)$, $E=L^{\, p}(G)$ and $J=\B(A, E)$. By Proposition \ref{coretract} there exists $R \in \AB(\widetilde{J}, E)$ with $R\circ \widetilde{\Pi}=I_{E}$. For a compact, non-null set $V\subset G$ we define a linear functional $\Lambda_{V}$ on $A'$ by
$$
\left\langle \lambda, \Lambda_{V}\right\rangle=\frac{1}{m(V)}\int_{V}R(\chi_{V}\otimes \lambda)(t) \dd m(t) \quad (\lambda \in L^{\, \infty}(G))\,.
$$
For $\lambda \in A'$ we have
$$
\mod{\left\langle \lambda, \Lambda_{V}\right\rangle}\leq \norm{R(\chi_{V}\otimes \lambda)}_{p}\norm{\chi_{V}/m(V)}_{p'}\leq C\norm{\lambda}_{\infty}\norm{\chi_{V}}_{p}\norm{\chi_{V}/m(V)}_{p'}=C\norm{\lambda}_{\infty}\,,
$$
and so $\Lambda_{V} \in A''$ with $\norm{\Lambda_{V}}\leq C$. Let $\mathcal{F}$ be the family of compact, non-null neighbourhoods of $e$ in $G$, and set $V_{1}\leq V_{2}$ if $V_2 \subset V_1$. Then $(\mathcal{F}, \leq)$ is a directed set. Let $\Lambda$ be a weak-$*$ accumulation point in $A''$ of the bounded net $\{\Lambda_{V}: V \in \mathcal{F}\}$. We {\it claim} that $\Lambda$ is $(p,p)$-multi-invariant.

Clearly $\left\langle 1, \Lambda\right\rangle=1$ since for each $V \in \mathcal{F}$ we have
$$
\left\langle 1, \Lambda_{V}\right\rangle=\frac{1}{m(V)}\int_{V}R(\widetilde{\Pi}\chi_{V})(t) \dd m(t)=\frac{1}{m(V)}\int_{V}\dd m(t)=1\,.
$$

Take $n \in \N$ and a finite subset $\{s_{1},\ldots, s_{n}\} \subset G$. Let $U \in \B(A, E)$, and let $\mathbf{X}=(X_1,\ldots, X_n)$ be a measurable partition of $G$. Take $f_1,\ldots, f_n \in L^{\, p'}(G)_{[1]}$ with $\supp f_i \subset X_i\,\; (i \in \N_n)$. Choose $V \in \mathcal{F}$ such that the sets $s_{1}V, \ldots,s_{n}V$ are pairwise disjoint. Set
$$
T=\sum_{i=1}^{n}\chi_{s_{i}V}\otimes U'(f_i): A\rightarrow E\,.
$$
By Lemma \ref{5.3.5}, $T\in J$ and $\norm{T}\leq \norm{U}m(V)^{1/p}$.

For each $i \in \N_n$ we have
\begin{equation*}
\begin{split}
m(V)\la f_i, U''(s_{i}\cdot \Lambda_{V})\ra &=m(V)\left\langle U'(f_i), s_{i}\cdot \Lambda_{V}\right\rangle\\
&=\int_{V}R(\chi_{V}\otimes (U'(f_i)\cdot s_i))(t) \dd m(t) \\
&=\int_{V}R((s_i\cdot \chi_{V})\otimes U'(f_i))(s_{i}t) \dd m(t)\quad \text{(by ~\eqref{eq4.3})} \\
&=\int_{s_{i}V}R(\chi_{s_{i}V}\otimes U'(f_i))(t) \dd m(t) \\
&=\int_{s_{i}V}R(\chi_{s_{i}V}T)(t) \dd m(t)\,.
\end{split}
\end{equation*}
Hence by \holders inequality we have
$$
\mod{\left\langle U'(f_i), s_{i}\cdot \Lambda_{V}\right\rangle}\leq\norm{\chi_{s_{i}V}R(\chi_{s_{i}V}T)}_{p}m(V)^{1/p'-1}\,.
$$
Then by Proposition \ref{5.3.4} we have
\begin{equation*}
\begin{split}
\left(\sum_{i=1}^{n}\mod{\la f_i, U''(s_{i}\cdot \Lambda_{V})\ra}^{p}\right)^{1/p} &\leq \left(\sum_{i=1}^{n}\norm{\chi_{s_{i}V}R(\chi_{s_{i}V}T)}_{p}^{p}\right)^{1/p}m(V)^{1/p'-1} \\
&\leq C\norm{T}m(V)^{1/p'-1}\leq \norm{R}\norm{U}\,.
\end{split}
\end{equation*}
Therefore
$$
\lp\sum_{i=1}^{n}\mod{\la f_i, U''(s_{i}\cdot \Lambda)\ra}^{p}\rp^{1/p} =\lim_{V}\left(\sum_{i=1}^{n}\mod{\left\langle U'(f_i), s_{i}\cdot \Lambda_{V}\right\rangle}^{p}\right)^{1/p}\leq \norm{R}\,.
$$
Since this is true for all such collections $(f_i)$, we have
$$
\left(\sum_{i=1}^{n}\norm{\chi_{X_i}U''(s_{i}\cdot \Lambda)}_p^{p}\right)^{1/p}\leq \norm{R}\,.
$$
Since this is true for each measurable partition $\mathbf{X}$ and $U \in J_{[1]}$, by Lemma \ref{4.6.20}, we have
$$
\norm{(s_1\cdot \Lambda,\ldots, s_n\cdot \Lambda)}_n^{(p,p)}\leq \norm{R}\,.
$$
Therefore $C_{p,p}(\Lambda)\leq\norm{R}$, and $G$ is $(p,p)$-amenable.
\end{proof}

\subsection{The discrete case}
The proof of Theorem \ref{5.3.7} becomes much simpler when $G$ is discrete. Let $G$ be a group, and let $1<p<\infty$. We set $A=\ell^{\, 1}(G)$, $E=\ell^{\, p}(G)$ and $J=\B(A, E)$. We shall identify $J$ with a space of functions in $\C^{G\times G}$ via $$U(t,s)=U(\delta_{t})(s)\quad (s, t \in G,\, U \in J)\,.$$ With this identification we have
$$
(r*U)(s,t)=U(r^{-1}s, r^{-1}t)\quad (r, s, t \in G,\, U \in J)\,.
$$

The following is a special case of Proposition \ref{coretract}, but the proof becomes much more direct when $G$ is discrete.

\begin{proposition}\label{discreteCoretract}
Let $G$ be a group, and let $1<p<\infty$. Then $\ell^{\, p}(G)$ is injective in $\lmod$ if and only if the morphism $\widetilde{\Pi}:E\rightarrow \widetilde{J}$ is a coretraction.
\end{proposition}
\begin{proof}
The condition is necessary by \cite[VII.1.34]{HE93}. We shall show it is sufficient.

Suppose that there is a morphism $R \in \AB(\widetilde{J}, E)$ with $R\circ \widetilde{\Pi}=I_{E}$. Define $Q:J\rightarrow \widetilde{J}$ by
$$
Q(U)(a)=(a*U)(\delta_{e})\quad (a \in A,\, U \in J)\,.
$$
For $U \in J$ we have
$$
Q(U)(t, s)=(t*U)(e, s)=U(t^{-1}, t^{-1}s)\quad (t, s \in G)\,.
$$
Now for $r \in G$ we have
\begin{equation*}
\begin{split}
\left(r*Q(U)\right)(t, s)&=Q(U)(r^{-1}t, r^{-1}s)=U(t^{-1}r, t^{-1}s) \\
												&=\left(r\cdot U\right)(t^{-1}, t^{-1}s)=Q(r\cdot U)(t, s)\,,
\end{split}
\end{equation*}
and hence $Q$ is a left $A$-module morphism. For $x \in E$, we have
$$
Q(\Pi(x))(t, s)=\Pi(x)(t^{-1}, t^{-1}s)=(t^{-1}\cdot x)(t^{-1}s)=x(s)=\widetilde{\Pi}(x)(t, s)\quad (s, t \in G)\,.
$$
Hence $Q\circ \Pi=\widetilde{\Pi}$. Finally, we set $\rho=R\circ Q$. Then $\rho \in \AB(J, E)$ and $\rho \circ \Pi=I_{E}$, and so $E$ is injective in $\Amod$.
\end{proof}

\begin{proposition}\label{5.3.9} \ifmargin \marginpar[5.3.9]{5.3.9} \fi
Let $S$ be a set, and take $1 < p <\infty$. Let $R:\B(\ell^{\, 1}(S), \ell^{\,p}(S))\rightarrow \ell^{\, p}(S)$ be a bounded linear operator. Then
$$
\left(\sum_{s \in S}\left|R(\delta_{s}U)(s)\right|^{p}\right)^{1/p}\leq \norm{U}\norm{R}\quad (U \in \B(\ell^{\, 1}(S), \ell^{\,p}(S)))\,.
$$
\end{proposition}
\begin{proof}
This follows from Proposition \ref{5.3.4}.
\end{proof}

By equation \eqref{eq4.3} we have
\begin{equation}\label{eq4.10}
\delta_{e}\otimes (U'(\delta_s)\cdot s)=s^{-1}*[\delta_{s}\otimes U'(\delta_s)]=s^{-1}*(\delta_{s}U)\,.
\end{equation}

\begin{theorem}\label{5.3.10} \ifmargin \marginpar[5.3.10]{5.3.10} \fi
Let $G$ be a group, and take $p$ with $1<p< \infty$. Then $\ell^{\, p}(G)$ is injective in $\lmod$ if and only if $G$ is $(p,p)$-amenable.
\end{theorem}
\begin{proof}
We set $A=\ell^{\, 1}(G)$ and $E=\ell^{\, p}(G)$. Suppose first that $E$ is injective in $\Amod$. By Proposition \ref{discreteCoretract} there exists $R \in \AB(\widetilde{J}, E)$ with $R\circ \widetilde{\Pi}=I_{E}$. We define $\Lambda \in A''$ by
$$
\left\langle \lambda, \Lambda\right\rangle=R(\delta_{e}\otimes \lambda)(e)\quad (\lambda \in A')\,.
$$
We have
$$
\left\langle 1, \Lambda\right\rangle=R(\delta_{e}\otimes \chi_{G})(e)=R(\widetilde{\Pi}\delta_{e})(e)=\delta_{e}(e)=1\,.
$$

We {\it claim} that $\Lambda$ is a $(p,p)$-multi-invariant mean. Take $n \in \N$ and a finite subset $\{s_{1},\ldots, s_{n}\} \subset G$. Let $U \in \B(A, E)$, and let $\mathbf{X}=(X_1,\ldots, X_n)$ be a partition of $G$. Take $f_1,\ldots, f_n \in E'_{[1]}$ with $\supp f_i \subset X_i\,\; (i \in \N_n)$. Set
$$
T=\sum_{i=1}^{n}\delta_{s_i}\otimes U'(f_i):A\rightarrow E\,.
$$
By Lemma \ref{5.3.5} $T\in J$ and $\norm{T}\leq \norm{U}$.

For each $i \in \N_n$ we have
\begin{equation*}
\begin{split}
\la f_i, U''(s_{i}\cdot \Lambda)\ra&=\left\langle U'(f_i), s_{i}\cdot \Lambda\right\rangle=R(\delta_{e}\otimes (U'(f_i)\cdot s_i))(e)\\
&=R(\delta_{s_i}\otimes U'(f_i))(s_{i}) \quad \text{(by \eqref{eq4.10})} \\
&=R(\delta_{s_i}T)(s_i)\,.
\end{split}
\end{equation*}
Then by Proposition \ref{5.3.9} we have
$$
\lp\sum_{i=1}^{n}\mod{\la f_i, U''(s_{i}\cdot \Lambda)\ra}^{p}\rp^{1/p} =\left(\sum_{i=1}^{n}\mod{R(\delta_{s_i}T)(s_i)}^{p}\right)^{1/p}\leq \norm{R}\norm{U}
$$
Since this is true for all such collections $(f_i)$, we have
$$
\left(\sum_{i=1}^{n}\norm{\chi_{X_i}U''(s_{i}\cdot \Lambda)}_p^{p}\right)^{1/p}\leq \norm{R}\norm{U}\,.
$$
Since this is true for each partition $\mathbf{X}$ and $U \in J$, by Lemma \ref{4.6.20} we have $$\norm{(s_{1}\cdot\Lambda, \ldots, s_{n}\cdot\Lambda)}^{(p,p)}_{n}\leq \norm{R}\,.$$ Therefore $C_{p,p}(\Lambda)\leq \norm{R}$, and $G$ is $(p,p)$-amenable.

Conversely, suppose that $G$ is $(p,p)$-amenable, and let $\Lambda \in A''$ be a $(p,p)$-multi-invariant mean. For $U \in J$ define $R(U):G\rightarrow \C$ by
$$
R(U)(s)=\left\langle U'(\delta_s), s\cdot\Lambda\right\rangle=U''(s\cdot \Lambda)(s)\quad (s \in G)\,.
$$
It is easily checked that $R(U) \in E$ and $R \in \B(J, E)$ with $\norm{R}\leq C_{p,p}(\Lambda)$. For $r, s \in G$ we have
$$
R(r * U)(s)=U''(r^{-1}s\cdot \Lambda)(r^{-1}s)=\ls r\cdot R(U)\rs(s)\,.
$$
Therefore $R \in \AB(J, E)$. For $x \in E$ and $s \in G$ we have
$$
R(\widetilde{\Pi} x)(s)=\la 1, \Lambda\ra x(s)=x(s)\,.
$$
Therefore $R\circ \widetilde{\Pi}=I_{E}$. By Proposition \ref{discreteCoretract}, $E$ is injective in $\Amod$.
\end{proof}

Since $(p,p)$-amenability implies $(q,q)$-amenability for any $q\geq p$, we have the following corollary.

\begin{corollary}
Let $G$ be a group, and take $p$ with $1<p< \infty$. Suppose that $\ell^{\, p}(G)$ is injective in $\lmod$. Then $\ell^{\, q}(G)$ is injective in $\lmod$ for all $q\geq p$.\qed
\end{corollary}

\begin{remark}
There are natural quantitative versions of projectivity, injectivity and flatness. These were first explicitly introduced and studied in \cite{MCW}. Let $A$ be a Banach algebra, and let $E \in \Amod$ be injective. We set $\inj(E)=\inf \norm{\rho}$ where the infimum is taken over all right inverse morphisms $\rho$ to the canonical morphism $\Pi$.

Let $G$ be a $(p,p)$-amenable locally compact group, and set $C_{p,p}(G)=\inf C_{p,p}(\Lambda)$ where the infimum is taken over all $(p,p)$-multi-invariant means $\Lambda$. The number $C_{p,p}(G)$ is related to $\inj(L^p(G))$. A chase through the constants in this section shows that $C_{p,p}(G)\leq \inj(L^p(G))$, and $C_{p,p}(G)= \inj(L^p(G))$ if $G$ is discrete.

These constants become significant in light of a recent result of G. Racher. Racher has proved (by different methods to us) that, for a discrete group $G$, if $\ell^{\,2}(G)$ is injective with $\inj(\ell^{\, 2}(G))=1$, then $G$ must be amenable. It follows from this result and Theorem \ref{5.3.10} that a discrete group $G$ is $(2,2)$-amenable with $C_{2,2}(G)=1$ if and only if $G$ is amenable.

We {\it conjecture} that for any $p>1$, any $(p,p)$-amenable locally compact group is amenable.
\end{remark}

{\bf Acknowledgements}
This article is based on work from the author's PhD
thesis \cite{PRthesis} at the University of Leeds, which was supported by
an EPSRC grant. Particular thanks are due to H. G. Dales and M. E. Polyakov for sharing preprints and unpublished notes, and whose original ideas underlie this work.




\providecommand{\bysame}{\leavevmode\hbox to3em{\hrulefill}\thinspace}


\begin{thebibliography}{10}

\bibitem{COHN}
Donald~L.\ Cohn, \emph{Measure theory}, {Birkh\"auser}, 1980.

\bibitem{HGD}
H.~G. Dales, \emph{Banach algebras and automatic continuity}, London
  Mathematical Society Monographs {\bf 24} (Clarendon Press, Oxford), 2000.

\bibitem{DLS2}
H.~G. Dales, A.~T.-M. Lau, and D.~Strauss, \emph{The measure algebra and its
  second dual}, preprint.

\bibitem{DP08}
H.~G. Dales and M.~E. Polyakov, \emph{Multi-normed spaces}, preprint.

\bibitem{DP}
\bysame, \emph{Homological properties of modules over group algebras}, Proc.
  London Math. Soc. \textbf{(2) 63} (2004), 390--426.

\bibitem{DS}
N.~Dunford and J.~T. Schwartz, \emph{Linear operators, {P}art {I}. {G}eneral
  theory}, Interscience, New York and London, 1958.

\bibitem{HE86}
A.~Ya. Helemskii, \emph{The homology of {B}anach and topological algebras},
  Kluwer Academic Publishers, Dordrecht, 1986.

\bibitem{HE93}
\bysame, \emph{{B}anach and locally convex algebras}, Clarendon Press, Oxford,
  1993.

\bibitem{JAM}
G.~J.~O. Jameson, \emph{Summing and nuclear norms in {B}anach space theory},
  London Society Student Texts, Volume 8, Cambridge University Press, 1987.

\bibitem{JO72}
B.~E. Johnson, \emph{Cohomology in {B}anach algebras}, Memoirs American Math.
  Soc. \textbf{127} (1972).

\bibitem{LNR04}
A.~Lambert, M.~Neufang, and V.~Runde, \emph{Operator space structure and
  amenability for {\figa}-{T}alamanca--{H}erz algebras}, J. Functional Analysis
  \textbf{211 (1)} (2004), 245--269.

\bibitem{PA}
A.~L.~T. Paterson, \emph{Amenability}, Mathematical Surveys and Monographs,
  Volume. 29, American Math. Soc., Providence, Rhode Island, 1988.

\bibitem{PRthesis}
P.~Ramsden, \emph{Homological properties of semigroup algebras}, Thesis,
  University of Leeds, 2009.

\bibitem{RUDIN73}
W.~Rudin, \emph{Functional analysis}, McGraw-Hill, New York, 1973.

\bibitem{TA79}
M.~Takesaki, \emph{Theory of operator algebras {I}}, Springer Verlag, New York,
  1979.

\bibitem{MCW}
M.~C. White, \emph{Injective modules for uniform algebras}, Proc. London Math.
  Soc. \textbf{(3) 73} (1996), 155--184.

\end{thebibliography}
\end{document}